\documentclass{article}
\usepackage{amsfonts,amsmath,amsthm,amssymb,xspace}
\font\cursief=pzcmi

\newcommand{\sectie}[1]{\setcounter{equation}{0}\section{#1}}

\newcommand{\N}{\mathbb{N}}
\newcommand{\Z}{\mathbb{Z}}

\newcommand{\R}{\mathbb{R}}
\newcommand{\C}{\mathbb{C}}

\newcommand{\hoed}{\,\hat{\rule{0ex}{1ex}}\;}
\newcommand{\alhh}{\hat{\alpha}
\hspace{-.7ex}\hat{\rule{0ex}{1.45ex}}\hspace{.7ex}}

\newcommand{\deh}{\hat{\Delta}}
\newcommand{\vfih}{\hat{\vfi}}

\newcommand{\lah}{\hat{\la}}
\newcommand{\pih}{\lambda}
\newcommand{\sih}{\hat{\si}}
\newcommand{\sdeh}{\hat{\delta}}
\newcommand{\nab}{\nabla}
\newcommand{\Jh}{\hat{J}}

\newcommand{\nabh}{\hat{\nab}}
\newcommand{\cI}{{\cal I}}
\newcommand{\cL}{{\cal L}}

\newcommand{\cN}{{\cal N}}
\newcommand{\cM}{{\cal M}}

\newcommand{\cD}{{\cal D}}

\newcommand{\ot}{\otimes}
\newcommand{\la}{\Lambda}
\newcommand{\om}{\omega}
\newcommand{\io}{\iota}
\newcommand{\vfi}{\varphi}
\newcommand{\al}{\alpha}
\newcommand{\be}{\beta}
\newcommand{\Ga}{\Gamma}
\newcommand{\sde}{\delta}
\newcommand{\de}{\Delta}
\newcommand{\si}{\sigma}
\newcommand{\Mfi}{{\cal M}_{\vfi}}
\newcommand{\Nfi}{{\cal N}_{\vfi}}
\newcommand{\Npsi}{{\cal N}_{\psi}}
\newcommand{\Mpsi}{{\cal M}_\psi}

\newcommand{\cU}{{\cal U}}

\newcommand{\cT}{{\cal T}}
\newcommand{\cJ}{{\cal J}}

\newcommand{\recht}{\rightarrow}
\newcommand{\qu}{$(M,\de)$\xspace}
\newcommand{\deop}{\de \hspace{-.3ex}\raisebox{0.9ex}[0pt][0pt]{\scriptsize\fontshape{n}\selectfont op}}
\newcommand{\fix}{N^\al}
\newcommand{\Next}{N^+
\hspace{-1.5ex}\raisebox{-0.5ex}[0pt][0pt]{\scriptsize\fontshape{n}\selectfont ext}}
\newcommand{\Ext}{^+\hspace{-1.5ex}\raisebox{-0.5ex}[0pt][0pt]{\scriptsize\fontshape{n}\selectfont ext}}
\newcommand{\vna}{von Neumann algebra\xspace}
\newcommand{\Mh}{\hat{M}}
\newcommand{\cros}{M \, \mbox{$_\al$}\hspace{-.2ex}\mbox{$\ltimes$} \, N}
\newcommand{\dehop}{\deh \hspace{-.3ex}\raisebox{0.9ex}[0pt][0pt]{\scriptsize\fontshape{n}\selectfont op}}
\newcommand{\alh}{\hat{\alpha}}
\newcommand{\tekst}[1]{\quad\text{#1}\quad}
\newcommand{\Wtil}{\tilde{W}}
\newcommand{\nsf}{\mbox{n.s.f.}\xspace}
\newcommand{\te}{\theta}
\newcommand{\Mte}{\cM_\theta}
\newcommand{\Nte}{\cN_\theta}
\newcommand{\Si}{\Sigma}
\newcommand{\Op}{\raisebox{0.9ex}[0pt][0pt]{\scriptsize\fontshape{n}\selectfont op}\,}
\newcommand{\ga}{\gamma}
\newcommand{\tetil}{\tilde{\theta}}
\newcommand{\Nfih}{\cN_{\vfih}}
\newcommand{\Ntetil}{\cN_{\tetil}}
\newcommand{\late}{\Lambda_\theta}
\newcommand{\latil}{\tilde{\Lambda}}
\newcommand{\strong}{$\si$-strong$^*$\xspace}
\newcommand{\lspan}{\operatorname{span}}
\newcommand{\tetiln}{\tilde{\theta}_0}
\newcommand{\Ntetiln}{\cN_{\tetiln}}
\newcommand{\Jtil}{\tilde{J}}
\newcommand{\nabtil}{\tilde{\nabla}}
\newcommand{\sitil}{\tilde{\sigma}}
\newcommand{\Ttil}{\tilde{T}}
\newcommand{\Jte}{J_\theta}
\newcommand{\nabte}{\nabla_\theta}
\newcommand{\site}{\sigma^\theta}
\newcommand{\Mfih}{\cM_{\vfih}}
\newcommand{\Ad}{\operatorname{Ad}}
\newcommand{\trace}{{\operatorname{Tr}}}
\newcommand{\tetiltil}{\tilde{\theta}
\hspace{-.4ex}\tilde{\rule{0ex}{1.85ex}}\hspace{.4ex}}
\newcommand{\tildeke}{\,\tilde{\rule{0ex}{1ex}}\;}
\newcommand{\na}{\, \raisebox{-.5ex}[0pt][0pt]{\scriptsize $^\circ$} \,}
\newcommand{\kruisje}[1]{\, \mbox{$_{#1}$}\hspace{-.2ex}\mbox{$\ltimes$}
\,}
\newcommand{\Te}{\Theta}
\newcommand{\cV}{\hspace{-.2ex}\mbox{\cursief V}\hspace{.5ex}}
\newcommand{\Ker}{\operatorname{Ker}}
\newcommand{\dis}{\displaystyle}
\newcommand{\sla}{\lambda}
\newcommand{\dpr}{^{\prime\prime}}
\newcommand{\Ntil}{\tilde{N}}
\newcommand{\Vtil}{\tilde{V}}
\newcommand{\hogechi}{\raisebox{.45ex}[0pt][0pt]{$\chi$}}

\newtheorem{definition}{Definition}[section]
\newtheorem{proposition}[definition]{Proposition}
\newtheorem{lemma}[definition]{Lemma}
\newtheorem{theorem}[definition]{Theorem}
\newtheorem{corollary}[definition]{Corollary}
\newtheorem{counterex}[definition]{Counterexample}

\begin{document}
\begin{center}
\huge\bf The unitary implementation of a locally compact quantum group action
\end{center}

\bigskip

\begin{center}
Stefaan Vaes\footnote{Research Assistant of the
Fund for Scientific Research -- Flanders (Belgium)\ (F.W.O.)}

Department of Mathematics

KU Leuven

Celestijnenlaan 200B

B--3001 Leuven

Belgium

\medskip

e-mail : Stefaan.Vaes@wis.kuleuven.ac.be

\bigskip
{\bf May 2000}

\bigskip
\end{center}
\begin{abstract}
\noindent In this paper we study actions of locally compact quantum groups on
von Neumann algebras and prove that every action has a canonical unitary
implementation, paralleling Haagerup's classical result on the unitary
implementation of a locally compact group action. This result is an important
tool in the study of quantum groups in action. We will use
it in this paper to study subfactors and inclusions of von Neumann algebras.
When $\al$ is an action of the locally compact quantum group \qu on the von Neumann algebra $N$ we can
give necessary and sufficient conditions under which the inclusion $N^\al
\subset N \hookrightarrow \cros$ is a basic construction.
Here $N^\al$ denotes the fixed point algebra and $\cros$ is the crossed
product.
When $\al$ is an outer and integrable action on a factor $N$ we prove
that the inclusion $N^\al \subset N$ is irreducible, of depth~2
and regular, giving a converse to the results of Enock and Nest, \cite{EN} and \cite{E}.
Finally we prove the equivalence of minimal and outer actions and
we generalize the main theorem of Yamanouchi \cite{Ya} :
every integrable outer action with infinite fixed point
algebra is a dual action.
\end{abstract}
\section*{Introduction}
Building on the work of Kac \& Vainerman \cite{VK}, Enock \& Schwartz \cite{ESbook}, Baaj
\& Skandalis \cite{BS},
Woronowicz \cite{Wor5} and Van Daele \cite{VD}, a precise definition of
a locally compact quantum group was recently introduced by the authors in
\cite{KV2}, see \cite{KV1} and \cite{KVPnas} for an overview.
For an overview of the historical development of the theory we refer
to \cite{KVPnas} and the introduction of \cite{KV2}.
This theory provides a topological framework to study
quantum groups and it unifies locally compact groups, compact
quantum groups and Kac algebras.

Because classical groups are usually defined to act on a space it
is very natural to make a quantum group act on a quantum space,
which will be an algebra. In an algebraic framework the study of
Hopf algebras acting on algebras has been very useful.

On the other hand, actions of locally compact groups on von
Neumann algebras have always been an important topic in operator
algebra theory, see e.g. \cite{Con} and \cite{NT}. In these
works the importance of Haagerup's results on the canonical
implementation of locally compact group actions and his results on
the dual weight construction, cannot be overestimated. It is
simply used all the time, without noticing it. See \cite{H1} and
\cite{H2}.

Hence it seems natural to study more generally actions of locally
compact quantum groups on von Neumann algebras and to try to
develop the same machinery of canonical implementation and dual
weight construction. This is what is done in the first half of
this paper. We strongly believe that this will serve as an
important tool in several applications of locally compact quantum
groups. We already give some applications in the second half of
this paper. Other applications are given by Kustermans in \cite{K}. Both will be explained below.

The special case
of Kac algebra actions has been studied by Enock and Schwartz in
\cite{E0} and \cite{ES}. They obtained important results on
crossed products, with the biduality theorem as a major
achievement. But they never obtained a unitary implementation for
an arbitrary action and also Haagerup's theory of dual weights on
the crossed product could not be completely generalized. For
instance, it remained an open problem whether the crossed product
with a Kac algebra action on a von Neumann algebra is in standard
position on its natural Hilbert space. It should be mentioned that
in \cite{BS} also Baaj \& Skandalis obtain a biduality theorem for crossed
products with multiplicative unitaries.

A first attempt to obtain the unitary implementation of a Kac
algebra action was made by J.-L.~Sauvageot in \cite{Sau}.
Unfortunately his proof is wrong, and for this we refer to the
discussion in the beginning of section~\ref{sectiehoeveel}.

So in this paper we will define actions of a locally compact quantum
group on a von Neumann algebra and we will construct its unitary
implementation. We will also give a construction for the dual
weight on the crossed product and prove analogous results as those
about group actions obtained by Haagerup in \cite{H1} and
\cite{H2}. In particular we prove that the crossed product is in
standard position on its natural Hilbert space and we identify the
associated modular objects. Hence we do not only give a right
proof for the results of Sauvageot, but also we prove new results
on the dual weights which make them a workable and applicable tool,
and we work in the more general setting of locally
compact quantum groups.

In the second half of the paper we will give some applications of
these results in the theory of subfactors and inclusions of von
Neumann algebras. It has been proved by Enock and Nest in their
beautiful papers \cite{EN} and \cite{E} that every irreducible, depth~2 inclusion
of factors satisfying the regularity condition, can be obtained as
$N^\al \subset N$ where $\al$ is an outer action of a locally
compact quantum group on the factor $N$ and $N^\al$ is the fixed
point algebra. We show in this paper that the action $\al$ is
always integrable and that, conversely, for every outer and
integrable action $\al$ on a factor $N$ the inclusion $N^\al
\subset N$ is irreducible, of depth~2 and regular. So we obtain a
converse to the results of Enock and Nest. The same result is stated for the special case of a dual
Kac algebra action in \cite[11.14]{EN}, but not proved.
While doing this, we study more generally the problem when the
inclusion $N^\al \subset N \hookrightarrow \cros$ is a basic
construction, and here $\cros$ denotes the crossed product.

As a final application of our results we prove the equivalence of
outerness and minimality of a locally compact quantum group action,
under the integrability condition. We also generalize the main
theorem of Yamanouchi \cite{Ya} to actions of arbitrary locally
compact quantum groups: when working on separable Hilbert spaces,
we prove that every integrable outer action with infinite fixed
point algebra is a dual action.

It should also be mentioned that our results on the unitary
implementation of a locally compact quantum group action are
already applied in a recent paper by Kustermans (see \cite{K})
in which he constructs induced corepresentations of locally compact
quantum groups. Taking into account the importance of induced
representations of locally compact groups, it is clear that the
results of Kustermans serve as a major motivation for our work.

{\bf Acknowledgment : }I would like to thank both J.~Kustermans
and L.~Vainerman for their valuable advice, interesting and
stimulating discussions and useful suggestions on this topic of
locally compact quantum groups in action.

\section*{Definitions and notations}
The whole of this paper will rely heavily on the modular theory of
von Neumann algebras. Throughout the text we will not make efforts
to give references to the original papers, but we will use
\cite{S} as a general reference.

When $\te$ is a normal, semifinite and faithful (we say \nsf)
weight on a von Neumann algebra $N$, one can make the so-called
GNS-construction $(K_\te,\pi_\te,\la_\te)$, where $K_\te$ is a
Hilbert space, $\pi_\te$ is a normal representation of $N$ on
$K_\te$ and $\la_\te: \Nte \recht K_\te$ is a linear map
satisfying $\pi_\te(x) \la_\te(y)= \la_\te(xy)$ for all $x \in N$
and $y \in \Nte$. Further $\la_\te(\Nte)$ is dense in $K_\te$.
Here $\Nte$ is the left ideal in $N$ defined by $\{x \in N \mid
\te(x^* x) < \infty \}$. The representation $\pi_\te$ is faithful
and often we will identify $N$ and $\pi_\te(N)$. Then we will use
the expression: let us represent $N$ on the GNS-space of $\te$
such that $(K_\te,\io,\la_\te)$ is a GNS-construction.

We will use several standard notations and results of modular
theory. We write
$$\Mte^+=\{x \in N^+ \mid \te(x) < \infty \}$$
and we denote with $(\site_t)_t$ the modular automorphism group of
$\te$. Further we denote with $\cT_\te$ the Tomita algebra defined
by
$$\cT_\te = \{x \in N \mid x \tekst{is analytic with respect to}
(\site) \tekst{and} \site_z(x) \in \Nte \cap \Nte^* \tekst{for
all} z \in \C \}.$$
Given a GNS-construction $(K_\te,\pi_\te,\la_\te)$ we define as
usual the modular conjugation $\Jte$ and the modular operator
$\nabte$. Recall that
$$\Jte \nabte^{1/2} \late(x)=\late(x^*)$$
for all $x \in \Nte \cap \Nte^*$ and $\late(\Nte \cap \Nte^*)$ is
a core for $\nabte^{1/2}$.

When $\te_1$ and $\te_2$ are \nsf weights on $N$ we denote with
$([D \te_1 : D \te_2]_t)_{t \in \R}$ the Connes cocycle as defined
in e.g. \cite[3.1]{S}.

Often we will make use of operator valued weights. When $N$ is a
von Neumann algebra we denote with $\Next$ the extended positive
part of $N$ as introduced by Haagerup in \cite{H4}, see e.g. \cite[11.1]{S}.
For the notion of operator valued weights we refer to
\cite{H4} or \cite[11.5]{S}. We will denote with $\langle
\cdot , \cdot \rangle$ the composition of elements of $\Next$ and
$N_*^+$. When $T$ is an operator valued weight we denote with
$\cN_T$ the left ideal of elements $x$ such that $T(x^*x)$ is
bounded.

All tensor products in this paper are either von Neumann algebra
tensor products or tensor products of Hilbert spaces. This will
always be clear from the context. We will use $\si$ to denote the
flip map on a tensor product $A \ot B$ of von Neumann algebras and
$\Si$ to denote the flip map on a tensor product $H \ot K$ of
Hilbert spaces. When $K$ is a Hilbert space and $\xi \in K$ we
denote with $\te_\xi$ the operator in $B(\C,K)$ given by
$\te(\sla) = \sla \xi$ for all $\sla \in \C$. When $H$ is a
Hilbert space and $\xi,\eta \in H$ we denote with $\om_{\xi,\eta}$
the usual vector functional in $B(H)_*$ given by
$\om_{\xi,\eta}(x) = \langle x \xi, \eta \rangle$. We use
$\om_\xi$ as a shorter notation for $\om_{\xi,\xi}$. We will
denote with $\cD(T)$ the domain of a (usually densily defined) map
$T$.

Throughout this paper the pair \qu will denote a (von Neumann
algebraic) locally compact quantum group. This means that
\begin{itemize}
\item $M$ is a von Neumann algebra and $\de : M \recht M \ot M$ is
a normal and unital $*$-homomorphism satisfying coassociativity:
$(\de \ot \io)\de = (\io \ot \de)\de$.
\item There exist \nsf weights $\vfi$ and $\psi$ on $M$ such that
\begin{itemize}
\item $\vfi$ is left invariant in the sense that $\vfi \bigl( (\om \ot
\io)\de(x) \bigr) = \vfi(x) \om(1)$ for all $x \in \Mfi^+$ and $\om \in
M_*^+$.
\item $\psi$ is right invariant in the sense that $\psi \bigl( (\io \ot
\om)\de(x) \bigr) = \psi(x) \om(1)$ for all $x \in \Mpsi^+$ and $\om
\in M_*^+$.
\end{itemize}
\end{itemize}
We refer to \cite{KV2}, \cite{KV3} and \cite{KV1} for the theory
of locally compact quantum groups in either the von Neumann
algebra or C$^*$-algebra language. It should be stressed that in
\cite{KV2} there is given a definition of a locally compact
quantum group in the C$^*$-algebra framework, and it is proven
that one can associate with this a locally compact quantum group
in the von Neumann algebra framework. In \cite{KV3} the above
definition of a von Neumann algebraic locally compact quantum
group is given and it is shown how to associate with it a
C$^*$-algebraic locally compact quantum group.

One can then prove the existence of a \strong closed map $S$ on
$M$, called the antipode, such that for all $a,b \in \Nfi$ we have
$$(\io \ot \vfi) \bigl( \de(a^*)(1 \ot b) \bigr) \in \cD(S) \tekst{and} S \bigl( (\io \ot \vfi) \bigl( \de(a^*)(1 \ot
b) \bigr) \bigr) = (\io \ot \vfi)\bigl( (1 \ot a^*) \de(b) \bigr).$$
Moreover, such elements $(\io \ot \vfi) \bigl( \de(a^*)(1 \ot b) \bigr)$ span a
\strong core for $S$. Then there exists a unique one-parameter
group $(\tau_t)_{t \in \R}$ of automorphisms of $M$ and a unique
$*$-anti-automorphism $R$ of $M$ such that
$$S = R \tau_{-i/2} \quad\quad R^2 = \io \quad \quad R \tau_t =
\tau_t R \quad\text{for all}\quad t \in \R.$$
We call $\tau$ the scaling group of \qu and $R$ the unitary
antipode. One refers to the expression $S = R \tau_{-i/2}$ as the
polar decomposition of the antipode.

Next one can prove that $\de R = \si(R \ot R) \de$, where $\si$
denotes the flip map on $M \ot M$. One can also prove that the
left and right invariant weights on \qu are unique up to a
positive scalar. So $\psi$ and $\vfi R$ are proportional and we
can suppose from the beginning that $\psi = \vfi R$. We denote
with $(\si_t)_{t \in \R}$ the modular group of $\vfi$. Then there
exists a unique self-adjoint, strictly positive operator $\sde$
affiliated with $M$ and satisfying $\si_t(\sde) = \nu^t \sde$ and
$\psi = \vfi_\sde$, where $\nu > 0$ is a real number. Formally
this means that $\psi(x) = \vfi(\sde^{1/2}x \sde^{1/2})$ and for
an exact definition of $\vfi_\sde$ we refer to \cite[1.5]{SV}. We call $\sde$ the
modular element of \qu. The number $\nu > 0$ is called the scaling
constant of \qu.

Let us represent $M$ on the GNS-space $H$ of $\vfi$ such that
$(H,\io,\la)$ is a GNS-construction for $\vfi$. Let $(H,\io,\Ga)$
be the canonical GNS-construction for $\psi=\vfi_\sde$ as defined
in \cite[7.2]{KV2}. Then one can define unitaries $W$ and $V$ in
$B(H \ot H)$ such that
\begin{align*}
W^* (\la(a) \ot \la(b)) &= (\la \ot \la) \bigl( \de(b)(a \ot 1) \bigr)
\tekst{for all} a,b \in \Nfi \\
V(\Ga(a) \ot \Ga(b)) &= (\Ga \ot \Ga) \bigl( \de(a)(1 \ot b) \bigr) \tekst{for
all} a,b \in \Npsi.
\end{align*}
Here $\la \ot \la$ and $\Ga \ot \Ga$ denote the canonical GNS-maps
for the tensor product weights $\vfi \ot \vfi$ and $\psi \ot
\psi$. Then $W$ and $V$ are multiplicative unitaries, which means
that they satisfy the pentagonal equation
$$W_{12} W_{13} W_{23} = W_{23} W_{12}.$$
Denoting with $^-$ the \strong closure we have that
$$M = \{(\io \ot \om)(W) \mid \om \in B(H)_* \}^- \tekst{and}
\de(x) = W^* (1 \ot x) W \tekst{for all} x \in M.$$

We will denote with $J$ and $\nabla$ the modular conjugation and
modular operator of $\vfi$ in the GNS-construction $(H,\io,\la)$.

Finally we describe how to define the dual locally compact quantum
group $(\Mh,\deh)$. Define the von Neumann algebra $\Mh$ as
follows, where again $^-$ denotes the \strong closure.
$$\Mh = \{(\om \ot \io)(W) \mid \om \in M_* \}^-.$$
Then one can define a comultiplication $\deh$ on $\Mh$ by
$$\deh(y) = \Si W (y \ot 1) W^* \Si \tekst{for all} y \in \Mh$$
where $\Si$ denotes the flip map on $H \ot H$. When $\om \in M_*$
we define $\pih(\om) = (\om \ot \io)(W)$. Of course $M_*$ should be thought of
as the L$^1$-functions on the quantum group \qu, and then $\lambda$ is the left
regular representation.
We also define
$$\cI = \{\om \in M_* \mid \;\text{there exists}\quad \eta \in H
\tekst{such that} \om(x^*)= \langle \eta,\la(x) \rangle \tekst{for
all} x \in \Nfi \}.$$
Such a $\eta \in H$ is necessarily uniquely determined and will be
denoted with $\xi(\om)$. Then there exists a unique \nsf weight
$\vfih$ on $\Mh$ with GNS-construction $(H,\io,\lah)$ such that
$\pih(\cI) \subset \Nfih$, $\pih(\cI)$ is a \strong--norm core for
$\lah$ and $\lah(\pih(\om)) = \xi(\om)$ for all $\om \in \cI$.

Then $(\Mh,\deh)$ will be a locally compact quantum group, and
having fixed the GNS-construction $(H,\io,\lah)$ for $\vfih$ we
can now repeat the story about $(M,\de)$ and obtain $(\sih_t)_{t
\in \R}$,$\sdeh$,$\hat{W}$,$\hat{V}$,$\Jh$ and
$\hat{\nabla}$. For all kinds of formulas relating these objects
we refer to \cite{KV3}. We only mention that
\begin{align*}
\Jh J &= \nu^{i/4} J \Jh \\
\nabh^{it} x \nabh^{-it} = \tau_t(x) \quad&\text{and}\quad \Jh x^*
\Jh = R(x) \tekst{for all} x \in M, t \in \R.
\end{align*}

Finally we denote with $(M,\de)\Op$ the opposite locally compact
quantum group $(M,\deop)$ where $\deop = \si \de$. Further we
define $(M,\de)'=(M',\de')$ where
$$\de'(x)=(J \ot J) \de(JxJ) (J \ot J)$$
for all $x \in M'$ and we call $(M,\de)'$ the commutant locally compact quantum group.
Then one can prove that $$(M,\de)\Op\hoed = (M,\de)\hoed ^\prime$$
and for this we refer to \cite{KV3}.

\sectie{Actions of locally compact quantum groups}
In this section we define actions of locally compact quantum
groups on von Neumann algebras and we construct an important tool:
the canonical operator valued weight from the von Neumann algebra
on which we act to the fixed point algebra, obtained by
integrating out the action.

\begin{definition} \label{11}
Let $N$ be a von Neumann algebra. A normal, injective and unital
$*$-homomorphism $\al:N \recht M \ot N$ will be called a left action of \qu on
$N$ when $(\io \ot \al)\al = (\de \ot \io)\al$.

A normal, injective and unital $*$-homomorphism $\al:N \recht N \ot M$ will be
called a right action of \qu on $N$ when $(\al \ot \io)\al = (\io \ot \de)\al$.
\end{definition}
In this paper we will only work with left actions and so we drop the predicate
\emph{left}. When $\al$ is a right action, $\si \al$ will be a left action of
$(M,\deop)$ on $N$, where $\si$ denotes the flip map from $N \ot M$ to $M \ot
N$ and $\deop$ denotes the opposite comultiplication. So it is not a real
restriction to work only with left actions. It should be observed
that in \cite{E0} and \cite{ES} they work with right actions.

\begin{definition}\label{12}
When $\al : N \recht M \ot N$ is an action of \qu on $N$ we define the
fixed point algebra $\fix$ as
$$\fix=\{ x \in N \mid \al(x)=1 \ot x\}.$$
\end{definition}
It is clear that $\fix$ is a von Neumann subalgebra of $N$.

Recall that $\Next$ denotes the extended positive part of $N$.

\begin{proposition}\label{13}
Let $N$ be a \vna and $\al$ an action of \qu on $N$. For every $x \in N^+$ the
element $T_\al(x) = (\psi \ot \io)\al(x)$ of $\Next$ belongs to $(N^\al)\Ext$.
Further $T_\al$ is a normal, faithful operator valued weight from $N$ to
$N^\al$.
\end{proposition}
\begin{proof}
Let $x \in N^+$ and $\om \in (M \ot N)^+_*$. Denote with $\langle \cdot,\cdot
\rangle$ the composition of an element in $\Next$ and an element in $N^+_*$.
Then by definition of the operator valued weight $\psi \ot \io$ we
get
\begin{align*}
\langle T_\al(x), \om \al \rangle = \langle (\psi \ot \io)\al(x), \om \al \rangle &= \psi \bigl( (\io \ot \om\al) \al(x) \bigr) \\
&=\psi \bigl( (\io \ot \om)(\de \ot \io)\al(x) \bigr) \\
&=\langle (\psi \ot \io \ot \io) \bigl( (\de \ot \io) \al(x) \bigr) , \om \rangle.
\end{align*}
By the right invariant version of \cite[3.1]{KV3} we get that
$$\langle T_\al(x), \om \al \rangle = \langle 1 \ot (\psi \ot \io)\al(x),\om
\rangle = \langle 1 \ot T_\al(x),\om \rangle.$$
From this it follows that $\al(T_\al(x)) = 1 \ot T_\al(x)$. So we get that
$T_\al(x) \in (N^\al)\Ext$.

If $x \in N^+$ and $a \in N^\al$ we have
$$\langle T_\al(axa^*),\om \rangle = \langle (\psi \ot \io) \bigl( (1 \ot a)\al(x)(1
\ot a^*) \bigr),\om \rangle = \psi \bigl( (\io \ot a^* \om a)\al(x) \bigr) = \langle
T_\al(x),a^* \om a \rangle.$$
So we get that $T_\al$ is indeed an operator valued weight. Because both $\al$
and $\psi \ot \io$ are faithful and normal, also $T_\al$ is faithful and normal.
\end{proof}

One should observe that the same result is stated and used in \cite{ES} for Kac
algebra actions. Their proof contains a gap because they do not have an
invariance formula like \cite[3.1]{KV3}, which is indispensable. For Kac algebra actions
this was repaired by Zsid\'o (see \cite{Zs}, see also \cite[18.18 and 18.23]{S}). Also in the case of a group action this was a non-trivial
problem (see \cite{H2}). The simpler proof of Zsid\'o for this invariance formula
only works in the Kac algebra setting, where the scaling group is trivial. I would like to thank prof. Enock for
the discussion on this topic.

\begin{definition} \label{14}
An action $\al$ of \qu on a \vna $N$ is called integrable when the operator
valued weight $T_\al$ defined in proposition~\ref{13} is semifinite.
\end{definition}

We will now introduce the well known concept of cocycle equivalent actions
(cfr. \cite[I.6]{E0}).
\begin{definition} \label{15}
Let $\al$ be an action of \qu on the \vna $N$. A unitary $U \in M \ot N$ is
called an $\al$-cocycle if
$$(\de \ot \io)(U)=U_{23} (\io \ot \al)(U).$$
It is clear that in this case the formula
$$\beta(x)=U\al(x) U^* \quad\text{for all}\quad x \in N$$
defines an action $\be$ of \qu on $N$.

Two actions $\al$ and $\be$ of \qu on $N$ are called cocycle equivalent if
there exists an $\al$-cocycle $U$ such that $\be$ is given by the formula
above.
\end{definition}

It is easy to see that $U^*$ is a $\be$-cocycle when $U$ is an $\al$-cocycle and when $\be(x)=U\al(x) U^*$ for all $x \in
N$.

\sectie{Crossed products, the dual action and the duality theorem}
In this section we fix an action $\al$ of a locally compact quantum group \qu on a \vna
$N$. We will define the crossed product $\cros$ in a similar way
as in \cite{E0}. We will also state some classical theorems
concerning crossed products, the biduality theorem being the major one,
but we will omit proofs because they
are completely analogous to the proofs of \cite{E0} and \cite{ES}.
See also \cite{E2}, where the results of \cite{E0} and \cite{ES}
are generalized to actions of Woronowicz algebras.
\begin{definition}\label{21}
We define the crossed product of $N$ with respect to the action $\al$ of \qu on $N$ as the von Neumann subalgebra of
$B(H) \ot N$ generated by $\al(N)$ and $\Mh \ot \C$. We denote this crossed product with $\cros$. So we have
$$\cros = (\al(N) \cup \Mh \ot \C)^{\prime\prime}.$$
\end{definition}
We will now define in the usual way the dual action, which will be an action of $(\Mh,\dehop)$ on $\cros$.
\begin{proposition} \label{22}
There exists a unique action $\alh$ of $(\Mh,\dehop)$ on $\cros$ such that
\begin{align*}
\alh(\al(x)) & = 1 \ot \al(x) \tekst{for all} x \in N \\
\alh(a \ot 1)&=\dehop(a) \ot 1 \tekst{for all} a \in \Mh.
\end{align*}
Moreover when we denote with $\Wtil$ the unitary $(J \ot J)\Si W \Si (J \ot J)$ we have
$$\alh(z) = (\Wtil \ot 1)(1 \ot z)(\Wtil^* \ot 1) \tekst{for all} z \in \cros.$$
\end{proposition}
As we already mentioned, Enock and Schwartz deal with right actions in \cite{E0}
and \cite{ES}. Hence they also give a slightly different definition
for the crossed product, but in fact our definition agrees with
theirs.
When $\al$ is a right action of \qu on $N$ they define $N \rtimes_\al M$ to be $(\al(N) \cup \C \ot
\Mh^{\prime})^{\prime \prime}$. This is in accordance with our definition, because $\si \al$ is a left action of
$(M,\deop)$ on $N$. The dual of $(M,\deop)$ is $(\Mh',\deh')$, which gives
$$M \kruisje{\si \al} N = (\si \al(N) \cup \Mh' \ot \C)^{\prime \prime}.$$
So we have $N \rtimes_\al M = \si(M \kruisje{\si \al} N)$, which
shows that both definitions in fact agree.

Let us introduce the following concept, which will be needed later on. See also \cite[III.1]{ES}.
\begin{definition}\label{23}
Let $\rho$ be a self-adjoint, strictly positive operator
affiliated with $M$. Then
a \nsf weight $\te$ on $N$ is called $\rho$-invariant if
$$\te \bigl( (\om_{\xi,\xi} \ot \io)\al(x) \bigr)=\|\rho^{1/2} \xi \|^2 \te(x)$$
for all $x \in \Mte^+$ and $\xi \in \cD(\rho^{1/2})$.
\end{definition}
We will always work with $\sde^{-1}$-invariant weights, where
$\sde$ is the modular element of the locally compact quantum group
in action.

Then the following result can be proved as in \cite[2.9]{E2}. For
the last statement of the next proposition observe that
$\tau_t(\sde)=\sde$ and so the self-adjoint operators $\sde$ and
$\nabh$ commute strongly. Hence their product $\sde \nabh$ is
closable.
\begin{proposition} \label{24}
When $\te$ is a \nsf $\sde^{-1}$-invariant weight on $N$ with
GNS-construction $(H_\te,\pi_\te,\la_\te)$, then there exists a unique
unitary $V_\te \in M \ot B(H_\te)$ such that for all $\xi \in
\cD(\sde^{1/2})$, $\eta \in H$ and $x \in \Nte$
$$(\om_{\xi,\eta} \ot \io)(V_\te) \la_\te(x)=
\la_\te \bigl( (\om_{\sde^{1/2} \xi,\eta} \ot \io)\al(x) \bigr).$$
Denote with $\Jte$ and $\nabte$ the modular conjugation and
modular operator of $\te$.
Then $V_\te$ satisfies
\begin{align*}
(\de \ot \io)(V_\te) &=V_{\te\, 23}V_{\te\, 13} \\
(\io \ot \pi_\te)\al(x) &= V_\te (1 \ot \pi_\te(x))V_\te^*
\tekst{for all} x \in N \\
V_\te(\Jh \ot J_\te) &= (\Jh \ot J_\te) V_\te^* \\
V_\te(Q \ot \nab_\te) &=(Q \ot \nab_\te) V_\te \tekst{where}
Q \;\text{ is the closure of }\; \sde \nabh.
\end{align*}
\end{proposition}

The following result is crucial (see \cite[2.8]{E2}).
\begin{proposition} \label{25}
\begin{itemize}
\item Let $\al$ be an integrable action of \qu on $N$ and denote with
$T_\al$ the operator valued weight defined in
proposition~\ref{13}. Let $\mu$ be a
\nsf weight on $N^\al$. Then $\mu \na T_\al$ is a
$\sde^{-1}$-invariant weight on $N$.
\item Every dual action is integrable.
\end{itemize}
\end{proposition}
With these results at hand one can copy the proofs of \cite{ES} to obtain
the well known biduality theorem.
Before we state this theorem we have to clarify some terminology.
The dual action $\alh$ is an action of $(M,\de)\hoed\Op$ on
$\cros$. So we can make the double crossed product $\Mh
\kruisje{\alh} (\cros)$ in $B(H \ot H) \ot N$ and on this double
crossed product there is an action $\alhh$ of
$(M,\de)\hoed\Op\hoed\Op$. Now $(M,\de)\hoed\Op\hoed\Op =
(M,\de)^\prime\Op$ and we can define an isomorphism of locally
compact quantum groups
$$\cJ : (M,\de) \recht (M,\de)^\prime \Op$$
given by $\cJ(x)=\Jh J x J \Jh$ for all $x \in M$.

\begin{theorem}[Biduality theorem]\label{26}
\begin{enumerate}
\item We have $B(H) \ot N = \bigl( B(H) \ot \C \cup \al(N) \bigr)^{\prime\prime}$.
\item The map $\Phi$ from $B(H) \ot N$ to $B(H \ot H) \ot N$ defined by
$$\Phi(z)=(W \ot 1) (\io \ot \al)(z) (W^* \ot 1)$$
defines a $*$-isomorphism from $B(H) \ot N$ onto $\Mh \kruisje{\alh} (\cros)$,
satisfying
\begin{alignat*}{3}
\Phi(\al(x)) &=1 \ot \al(x) & &\tekst{for all} & x &\in N \\
\Phi(b \ot 1) &= \dehop(b) \ot 1 & &\tekst{for all} & b &\in \Mh \\
\Phi(y \ot 1) &= y \ot 1 \ot 1 & &\tekst{for all} & y &\in M'.
\end{alignat*}
In particular $\Phi(\cros) = \alh(\cros)$.
\item When we define
$$\mu = (\si \ot \io)(\io \ot \al) : B(H) \ot N \recht M \ot B(H) \ot N$$
then $\mu$ is an action of \qu on $B(H) \ot N$. The unitary $\Si V^* \Si \ot 1$
is a $\mu$-cocycle and the action $\ga$ of \qu on $B(H) \ot N$ defined by
$$\ga(z)= (\Si V^* \Si \ot 1) \mu(z) (\Si V \Si \ot 1) \tekst{for all} z \in B(H) \ot N$$
is isomorphic to the bidual action $\alhh$ of $(M,\de)\hoed\Op\hoed\Op$ on $\Mh \kruisje{\alh}
(\cros)$ in the following way:
$$\alhh(\Phi(z)) = (\cJ \ot \Phi)\ga(z) \tekst{for all} z \in B(H) \ot
N.$$
\end{enumerate}
\end{theorem}
With the help of the biduality theorem Enock and Schwartz were able to prove
the following crucial results, which remain true for actions of locally compact
quantum groups.
\begin{theorem}\label{27}
We have
\begin{align*}
(\cros)^{\alh} & = \al(N) \\
\al(N) &= \{z \in M \ot N \mid (\io \ot \al)(z) = (\de \ot \io)(z) \}.
\end{align*}
\end{theorem}

\sectie{The unitary implementation of a locally compact quantum group
action}\label{sectiehoeveel}
In this section we will define in a canonical way the unitary implementation of
a locally compact quantum group action. This will be a unitary
corepresentation of the quantum group, implementing the action and satisfying
some other properties. A same kind of result was obtained for Kac algebra
actions by Sauvageot in \cite{Sau}, but the proof of the fact that the implementation
is a corepresentation, is wrong. More precisely, Sauvageot's crucial
lemma~4.1 is false. I would like to thank prof. Sauvageot for the discussions on
this topic.

We will use a different technique to prove that the implementation is a
corepresentation. In the same time we will obtain some interesting results
concerning the dual weight on the crossed product $\cros$ given a weight on
$N$. We will also settle a problem which was left open in \cite{Sau}.

For integrable actions -- and in particular for dual actions -- we already
obtained an implementation in proposition~\ref{24}, as it was done by Enock and
Schwartz. Nevertheless it is desirable to have an implementation without the
integrability condition, first of all for reasons of elegance. But, more
importantly, one will need this general implementation result in several
applications. We refer to the introduction for a discussion.

Fix an action $\al$ of a locally compact quantum group \qu on a \vna $N$. In definition~\ref{21} and
proposition~\ref{22} we defined the crossed product $\cros$ and the dual action
$\alh:\cros \recht \Mh \ot (\cros)$. We already observed in
proposition~\ref{25} that $\alh$ is integrable. So we can define the \nsf operator
valued weight $T$ from $\cros$ to $(\cros)^{\alh}$ by
$$T(z) = (\vfih \ot \io \ot \io) \alh(z) \tekst{for all} z \in (\cros)^+.$$
For this, observe that $\alh$ is an action of $(\Mh,\dehop)$ and that $\vfih$ is the
right invariant weight on $(\Mh,\dehop)$. By theorem~\ref{27} we know that
$(\cros)^{\alh} = \al(N)$. So $T$ is an operator valued weight from $\cros$ to
$\al(N)$.

With this operator valued weight at hand, we can easily define the
dual weights on $\cros$. Nevertheless, to make dual weights a
workable tool, we need a concrete GNS-construction for them. The
structure of this section is then as follows. First we will
restrict the dual weight to a weight for which we can give a
GNS-construction (definition \ref{34}), then we use the restricted
weight to obtain the unitary implementation for the action
(definition \ref{36} and proposition \ref{37}) and finally we
prove that the restricted weight is in fact not a restriction, but
equal to the original dual weight (proposition \ref{39b}).
\begin{definition}\label{30}
Let $\al$ be an action of \qu on $N$. Denote with $T$ the \nsf operator
valued weight from $\cros$ to $\al(N)$ given by the formula above.
For every \nsf weight $\te$ on $N$, we define the dual
weight $\tetil$ on $\cros$ by the formula:
$$\tetil = \te \na \al^{-1} \na T.$$
\end{definition}

For the rest of this section we fix a \nsf weight $\te$ on $N$.
One can prove easily the following lemma.
\begin{lemma}\label{32}
For all $a \in \Nfih$ and $x \in N$ we have
$$\tetil \bigl( \al(x^*)(a^* a \ot 1) \al(x) \bigr) = \te(x^* x) \vfih(a^* a).$$
\end{lemma}
\begin{proof}
We have
$$\alh \bigl( \al(x^*)(a^* a \ot 1) \al(x) \bigr) = (1 \ot \al(x^*))(\dehop(a^*a) \ot 1) (1
\ot \al(x)).$$
Choose $\om \in (\cros)^+_*$. Define $\mu \in \Mh_*^+$ by
$\mu(b)=\om \bigl( \al(x^*)(b \ot 1) \al(x) \bigr)$ for all $b \in \Mh$. Then
$$\langle T \bigl( \al(x^*)(a^* a \ot 1) \al(x) \bigr),\om \rangle = \vfih \bigl( (\io \ot
\mu)\dehop(a^*a) \bigr) = \vfih(a^*a) \mu(1)=\vfih(a^*a) \om(\al(x^*x))$$
by invariance of $\vfih$. So we may conclude that
$$T \bigl( \al(x^*)(a^* a \ot 1) \al(x) \bigr)=\vfih(a^*a) \al(x^*x).$$
Then the result of the lemma follows immediately.
\end{proof}
From now on we will suppose that $N$ acts on the GNS-space of the \nsf weight
$\te$, such that $(K,\io,\late)$ is a GNS-construction for $\te$. We will
restrict the weight $\tetil$ in the sense of proposition~\ref{a3} of the appendix in order to obtain
a concrete GNS-construction. Fix a GNS-construction $(K_1,\pi_1,\la_1)$ for
$\tetil$. Because of the previous lemma we can define a unique isometry
$$\cV:H \ot K \recht K_1 \tekst{such that} \cV(\lah(a) \ot \late(x)) = \la_1 \bigl( (a \ot
1)\al(x) \bigr)$$
for all $a \in \Nfih$ and $x \in \Nte$.

Further we define
$$\cD_0=\lspan \{ (a \ot 1) \al(x) \mid a \in \Nfih, x \in \Nte \}.$$
Because we have the isometry $\cV$ at our disposal there is a well defined
linear map
$$\latil_0:\cD_0 \recht H \ot K : \latil_0 \bigl( (a \ot 1)\al(x) \bigr) = \lah(a) \ot
\late(x) \tekst{for all} a \in \Nfih,x \in \Nte.$$
Because $\la_1$ is \strong--norm closed, we can close $\latil_0$ for the
\strong--norm topology, and then we obtain a linear map $\latil:\cD \recht  H
\ot K$ satisfying $\cD \subset \Ntetil$ and $\cV \latil(z)=\la_1(z)$ for all $z
\in \cD$.

In order to apply proposition~\ref{a3}, we need the following lemma.
\begin{lemma} \label{33}
\begin{enumerate}
\item $\cD$ is a weakly dense left ideal in $\cros$.
\item For all $z \in \cros$ and $y \in \cD$ we have $\latil(zy)=z\latil(y)$.
\end{enumerate}
\end{lemma}
\begin{proof}
Choose $\xi \in H$ and $b \in \cT_\vfi$. Let $(e_i)_{i \in I}$ be an
orthonormal basis for $H$. Choose $x \in N$. Because $(\de \ot \io)\al(x)=(\io
\ot \al)\al(x)$ we have
$$(1 \ot \al(x))(W \ot 1) = (W \ot 1)(\io \ot \al)\al(x).$$
Hence applying $\om_{\xi,\la(b)} \ot \io \ot \io$ gives
$$\al(x)(\pih(\om_{\xi,\la(b)}) \ot 1) = \sum_{i \in I} (\pih(\om_{e_i,\la(b)})
\ot 1) \al \bigl( (\om_{\xi,e_i} \ot \io)\al(x) \bigr)$$
in the \strong topology.

Choose now $y \in \Nte$. For every finite subset $I_0 \subset I$ we have by
proposition~\ref{a1} that the element
\begin{align*}
z_{I_0} &:= \sum_{i \in I_0} (\pih(\om_{e_i,\la(b)})
\ot 1) \al \bigl( (\om_{\xi,e_i} \ot \io)\al(x)y \bigr) \\
\intertext{belongs to $\cD_0$ and}
\latil_0(z_{I_0}) &= \sum_{i \in I_0} \lah \bigl( \pih(\om_{e_i,\la(b)}) \bigr) \ot (\om_{\xi,e_i} \ot
\io)\al(x)\late(y) \\
&=\sum_{i \in I_0} J \si_{i/2}(b)J e_i \ot (\om_{\xi,e_i} \ot \io)\al(x)
\late(y) \\
&= (J \si_{i/2}(b)J \ot 1)(P_{I_0} \ot 1) \al(x) (\xi \ot \late(y))
\end{align*}
where $P_{I_0}$ is the projection on $\lspan\{e_i\mid i \in I_0\}$. So we get
that the net $(z_{I_0})$ converges \strong to the element
\begin{align*}
z &:=\al(x)(\pih(\om_{\xi,\la(b)}) \ot
1)\al(y) \\
\intertext{and the net $(\latil(z_{I_0}))$ converges in norm to}
(J \si_{i/2}(b)J \ot 1) \al(x) (\xi \ot
\late(y)) &= \al(x) \bigl( J \si_{i/2}(b) J \xi \ot \late(y) \bigr) \\ &=
\al(x) \bigl( \lah(\pih(\om_{\xi,\la(b)})) \ot \late(y) \bigr) =
\al(x)\latil \bigl( (\pih(\om_{\xi,\la(b)}) \ot 1)\al(y) \bigr).
\end{align*}
Then we may conclude that $z \in \cD$ and
$$\latil(z)=\al(x)\latil \bigl( (\pih(\om_{\xi,\la(b)}) \ot 1)\al(y) \bigr).$$
Because the considered elements $\pih(\om_{\xi,\la(b)})$ form a \strong--norm
core for $\lah$ we conclude that for every $x \in N$ and $z \in \cD$ we have
$\al(x)z \in \cD$ and $\latil(\al(x)z) = \al(x) \latil(z)$.

It is easy to prove that for every $a \in \Mh$ and $z \in \cD$ we have $(a \ot
1)z \in \cD$ and $\latil \bigl( (a \ot 1)z \bigr) = (a \ot 1)\latil(z)$. From this follows
the lemma.
\end{proof}
We can now apply proposition~\ref{a3}.
\begin{definition} \label{34}
There is a unique \nsf weight $\tetiln$ on $\cros$ such that $\Ntetiln = \cD$
and such that $(H \ot K, \io, \latil)$ is a GNS-construction for $\tetiln$.
\end{definition}
Later on we will prove that in fact $\tetiln=\tetil$. This question was left
open in the Kac algebra case considered by Sauvageot. In
applications the equality $\tetiln=\tetil$ is indispensable,
e.g. proposition \ref{54} cannot be proved without knowing the
GNS-construction of $\tetil$, which amounts to the equality
$\tetiln=\tetil$.

Let us fix some modular notations.
\begin{definition} \label{35}
We denote with $\Jtil$ and $\nabtil$ the modular conjugation and modular
operator of $\tetiln$ in the GNS-construction $(H \ot K,\io,\latil)$. We denote
with $\sitil$ the modular automorphism group of $\tetiln$ and we put
$\Ttil=\Jtil \nabtil^{1/2}$.

We denote with $\Jte$ and $\nabte$ the modular conjugation and modular
operator of $\te$ in the GNS-construction $(K,\io,\late)$, and with $\site$ the
modular automorphism group of $\te$.
\end{definition}
With this notations at hand we will now define the unitary implementation of
the action $\al$. Of course this terminology will only be justified after the
proofs of~\ref{37}, \ref{310} and \ref{44}.
\begin{definition} \label{36}
Define $U=\Jtil(\Jh \ot \Jte)$. Then $U$ is a unitary in $B(H \ot K)$ and it is
called the unitary implementation of $\al$.
\end{definition}
We will first prove the following result.
\begin{proposition}\label{37}
We have the following formulas:
\begin{enumerate}
\item $\al(x)=U(1 \ot x)U^* \tekst{for all} x \in N$.
\item $\sitil_t \na \al = \al \na \site_t \tekst{for all} t \in \R$.
\item $U(\Jh \ot \Jte)=(\Jh \ot \Jte)U^*$.
\end{enumerate}
\end{proposition}
Before we can prove this proposition we need the following lemma.
\begin{lemma}\label{39}
For all $y \in \cD(\site_{i/2})$ we have $\al(y) \in \cD(\sitil_{i/2})$ and
$$\Jtil \sitil_{i/2}(\al(y))^* \Jtil = 1 \ot \Jte \site_{i/2}(y)^* \Jte.$$
\end{lemma}
\begin{proof}
Choose $a \in \Nfih$ and $x \in \Nte$. Then $xy \in \Nte$ and hence $(a \ot
1)\al(x)\al(y) \in \Ntetiln$ with
$$\latil \bigl( (a \ot 1) \al(x)\al(y) \bigr) = \lah(a) \ot \late(xy) =
(1 \ot \Jte \site_{i/2}(y)^* \Jte)\latil \bigl( (a \ot 1) \al(x) \bigr).$$
Because $\cD_0$ is a \strong--norm core for $\latil$ we may conclude that for
every $z \in \Ntetiln$ we have $z \al(y) \in \Ntetiln$ and
$$\latil(z \al(y)) = (1 \ot \Jte \site_{i/2}(y)^* \Jte)\latil(z).$$
Then the lemma follows immediately.
\end{proof}
\begin{proof}[{\bf Proof of proposition~\ref{37}.}]
Because $\site_{i/2}(y)^* = \site_{-i/2}(y^*)$ it follows from the previous
lemma that for every $y \in \cD(\site_{-i/2})$ we have $\al(y) \in
\cD(\sitil_{-i/2})$ and
$$\sitil_{-i/2}(\al(y)) = U(1 \ot \site_{-i/2}(y))U^*.$$
Taking the adjoint we may replace $-i/2$ by $i/2$ in the formula above.

Let now $y \in \cD(\site_{-i})$. Then we have $\al(y) \in \cD(\sitil_{-i/2})$
and
$$\sitil_{-i/2}(\al(y)) = U (1 \ot \site_{-i/2}(y))U^*.$$
Because $\site_{-i}(y) \in \cD(\site_{i/2})$ we also have $\al(\site_{-i}(y))
\in \cD(\sitil_{i/2})$ and
$$\sitil_{i/2} \bigl( \al(\site_{-i}(y)) \bigr) = U \bigl( 1 \ot \site_{i/2} ( \site_{-i}(y)) \bigr)U^* = U(1 \ot \site_{-i/2}(y))U^*.$$
So we get $\sitil_{-i/2}(\al(y))=\sitil_{i/2} \bigl( \al(\site_{-i}(y)) \bigr)$ and so
$\al(y) \in \cD(\sitil_{-i})$ with $\sitil_{-i}(\al(y)) =\al(\site_{-i}(y))$.
It now follows from the results of \cite[4.3 and 4.4]{H4} that $\sitil_t \na \al = \al \na \site_t$ for every $t \in \R$.

But then it follows that for all $y \in \cD(\site_{-i/2})$ we have
$\sitil_{-i/2}(\al(y)) = \al(\site_{-i/2}(y))$. Combining this with the formula
above we get
$$\al(\site_{-i/2}(y))=\sitil_{-i/2}(\al(y))=U(1 \ot \site_{-i/2}(y)) U^*.$$
By the density of such elements $\site_{-i/2}(y)$ we get that $\al(x)=U(1 \ot
x)U^*$ for all $x \in N$.

From the definition of $U$ follows immediately the final formula we had to
prove.
\end{proof}
Now we have gathered enough material to prove that $\tetiln=\tetil$. For this
we need the following lemma (cfr. \cite[VI.4]{ES}).
\begin{lemma}\label{39a}
Let $\al$ be an action of \qu on $N$. Let $\te_1$ and $\te_2$ be two
$\sde^{-1}$-invariant \nsf weights on $N$. Then $[D\te_2 : D \te_1]_t \in
N^\al$ for all $t \in \R$.
\end{lemma}
\begin{proof}
Denote with $M_2$ the \vna of $2 \times 2$-matrices over $\C$. Denote with
$e_{ij}$ the matrix units. Define
$$\ga:N \ot M_2 \recht M \ot N \ot M_2 : \ga=\al \ot \io.$$
Then $\ga$ is an action of \qu on $N \ot M_2$. Denote with $\te$ the balanced
weight on $N \ot M_2$ (see e.g. \cite[3.1]{S}) given by
$$\te \begin{pmatrix} x_{11} & x_{12} \\ x_{21} & x_{22} \end{pmatrix} = \te_1(x_{11}) +
\te_2(x_{22}).$$
It is immediately clear that $\te$ is $\sde^{-1}$-invariant for the action
$\ga$.

Let $t \in \R$. Denote with $\mu_t$ the automorphism of $M$ defined by $\mu_t = \si'_t \na
\si_{-t} \na \tau_t$. Here $(\si'_t)_{t \in \R}$ denotes the
modular automorphism group of $\psi$.
Then $\mu_t$ is implemented by $Q^{it} = \sde^{it}
\nabh^{it}$. It follows from proposition~\ref{24} that $\ga \na \site_t =
(\mu_t \ot \site_t) \na \ga$ for all $t \in \R$. In particular we have
\begin{align*}
\al([D \te_2 : D \te_1]_t) \ot e_{21} &= \ga( [D \te_2 : D \te_1]_t \ot
e_{21})=\ga \bigl( \site_t(1 \ot e_{21}) \bigr) \\
&=(\mu_t \ot \site_t)\ga(1 \ot e_{21}) = (\mu_t \ot \site_t)(1 \ot 1 \ot
e_{21}) = 1 \ot [D \te_2 : D \te_1]_t \ot e_{21}.
\end{align*}
So we get $[D \te_2 : D \te_1]_t \in N^\al$ for all $t \in \R$.
\end{proof}
Now we can prove the following interesting result. It is important
for technical reasons and we will need it in
section~\ref{sectie5}.
\begin{proposition}\label{39b}
Let $\te$ be a \nsf weight on $N$. Then the weights $\tetil$ and $\tetiln$
on $\cros$, defined in \ref{30} and \ref{34} are equal.
\end{proposition}
\begin{proof}
Recall that the dual action $\alh$ is an action of $(\Mh,\dehop)$ on $\cros$.
We claim that the weight $\tetiln$ is $\sdeh$-invariant. Observe that
$\sdeh^{-1}$ is the modular element of $(\Mh,\dehop)$ and that is the reason to
have $\sdeh$-invariance rather than $\sdeh^{-1}$-invariance.

To prove our claim, choose $a \in \Nfih$, $x \in \Nte$, $\xi \in
\cD(\sdeh^{1/2})$ and $\eta \in H$. Then define
$$z:=(\om_{\xi,\eta} \ot \io \ot \io)\alh \bigl( (a \ot 1) \al(x) \bigr) = \bigl( (\om_{\xi,\eta} \ot
\io)\dehop(a) \ot 1 \bigr) \al(x).$$
It follows from proposition~\ref{a2} of the appendix that
$$(\om_{\xi,\eta} \ot \io)\dehop(a) = (\io \ot \om_{\xi,\eta})\deh(a) \in \Nfih
\tekst{and} \lah \bigl( (\om_{\xi,\eta} \ot \io)\dehop(a) \bigr) = (\io \ot
\om_{\sdeh^{1/2}\xi,\eta})(\hat{V}) \lah(a).$$
So we may conclude that $z \in \Ntetiln$ and
$$\latil(z)=\bigl( (\io \ot \om_{\sdeh^{1/2}\xi,\eta})(\hat{V}) \ot 1 \bigr) \bigl( \lah(a) \ot
\late(x) \bigr) = \bigl( (\io \ot \om_{\sdeh^{1/2}\xi,\eta})(\hat{V}) \ot 1 \bigr) \latil \bigl( (a \ot
1)\al(x) \bigr).$$
Because $\cD_0$ is a \strong--norm core for $\latil$ we conclude that
$(\om_{\xi,\eta} \ot \io \ot \io)\alh(y) \in \Ntetiln$ for all $y \in \Ntetiln$
and
$$\latil \bigl( (\om_{\xi,\eta} \ot \io \ot \io)\alh(y) \bigr)= \bigl( (\io \ot \om_{\sdeh^{1/2}\xi,\eta})(\hat{V}) \ot
1 \bigr) \latil(y).$$
Because $\hat{V}$ is unitary we immediately get that $\tetiln$ is
$\sdeh$-invariant.

From proposition~\ref{25} it follows that $\tetil$ is $\sdeh$-invariant. Then
we conclude from lemma~\ref{39a} that $[D \tetiln: D\tetil]_t \in
(\cros)^{\alh}$ for all $t \in \R$. So by theorem~\ref{27} we can take
unitaries $u_t \in N$ such that $[D \tetiln: D\tetil]_t =\al(u_t)$ for all $t
\in \R$. From the theory of operator valued weights we know that
$\si_t^{\tetil} \na \al = \al \na \site_t$. Because $([D \tetiln:
D\tetil]_t)$ is a $\si^{\tetil}$-cocycle, we get that $(u_t)$ is a
$\site$-cocycle. By \cite[5.1]{S} we can take a (uniquely determined) \nsf weight $\rho$
on $N$ such that $[D\rho : D \te]_t = u_t$ for all $t \in \R$. With $\rho$ we
can define the \nsf weight $\tilde{\rho}$ on $\cros$ in the sense of definition \ref{30}. Then it follows from the
theory of operator valued weights that
$$[D\tilde{\rho} : D \tetil]_t = \al([D\rho:D\te]_t) = \al(u_t)=[D \tetiln:D
\tetil]_t$$
for all $t \in \R$. So $\tilde{\rho}=\tetiln$. Because $\tetiln$ is a
restriction of $\tetil$ we get that $\tilde{\rho}$ is a restriction of
$\tetil$.

Fix $a \in \Mfih^+$ with $\vfih(a)=1$. Choose $x \in \cN_\rho$. Then it follows
from lemma~\ref{32} that $\al(x^*) (a \ot 1) \al(x) \in \cM_{\tilde{\rho}}^+$
and
$$\tilde{\rho} \bigl( \al(x^*) (a \ot 1) \al(x) \bigr) = \rho(x^*x).$$
Because $\tilde{\rho}$ is a restriction of $\tetil$ we get that $\al(x^*) (a \ot 1)
\al(x)\in \cM_{\tetil}^+$ and
$$\tetil \bigl( \al(x^*) (a \ot 1) \al(x) \bigr) = \rho(x^*x).$$
Then it follows from lemma~\ref{32} that $\te(x^*x)=\rho(x^*x)$.
This means that $\rho$ is a restriction of $\te$.

Further we have, using the theory of operator valued weights in
the first equality and proposition \ref{37} in the last one,
$$\al \na \si^\rho_t = \si^{\tilde{\rho}}_t \na \al = \si^{\tetiln}_t \na
\al = \al \na \site_t.$$
So $\si^\rho_t = \si^\te_t$ for all $t \in \R$. Because $\rho$ is a
restriction of $\te$ we may conclude that $\rho=\te$ and then
$\tetil=\tilde{\rho}=\tetiln$.
\end{proof}
We want to conclude this section with the proof of the fact that
$U \in M \ot B(K)$. First we state the following lemma, which is
easily proved because $\latil$ is the closure of $\latil_0$.
Recall that $\Ttil = \Jtil \nabtil^{1/2}$.
\begin{lemma} \label{38}
Defining $\hat{T} = \Jh \nabh^{1/2}$, we have that the linear space
$$\lspan\{\al(x^*) (\eta \ot \late(y)) \mid x,y \in \Nte, \eta \in \cD(\hat{T}) \}$$
is a core for $\Ttil$ and
$$\Ttil \al(x^*)(\eta \ot \late(y)) = \al(y^*)(\hat{T}\eta \ot
\late(x))$$
for all $x,y \in \Nte$ and $\eta \in \cD(\hat{T})$.
\end{lemma}
\begin{proposition} \label{310}
We have $U \in M \ot B(K)$.
\end{proposition}
\begin{proof}
Let $t \in \R$. Because $\nabh^{it}$ implements the automorphism $\tau_t$ on $M$
we get that
$\Ad \nabh^{it}$ will also leave $M'$ invariant. So we can define
the automorphism group $(\mu_t)$ on $M$ by
$$\mu_t(x)=J\nabh^{it}JxJ\nabh^{-it} J \tekst{for all} x \in M, t \in \R.$$
So, for every $a \in \cD(\mu_{-i/2})$ we have $Ja J \nabh^{1/2}
\subset \nabh^{1/2} J \mu_{-i/2}(a)J$. Further we have
\begin{align*}
\mu_t(R(a)) &=J\nabh^{it}J\Jh a^* \Jh J \nabh^{-it} J = J
\nabh^{it} \Jh J a^* J \Jh \nabh^{-it} J \\
&= \Jh \mu_t(a^*) \Jh = R(\mu_t(a))
\end{align*}
for all $t \in \R$ and $a \in M$. Here we used the formula $\Jh J = \nu^{i/4} J \Jh$
stated in the beginning of the paper.

Let now $a \in \cD(\mu_{i/2})$, $x,y \in \Nte$ and $\eta \in
\cD(\hat{T})$, where $\hat{T} = \Jh \nabh^{1/2}$. Then
\begin{align*}
(JaJ \ot 1) \Ttil \al(x^*)(\eta \ot \late(y)) = (JaJ \ot 1)
\al(y^*) (\hat{T} \eta \ot \late(x)) \\
&= \al(y^*) \bigl( JaJ \Jh \nabh^{1/2} \eta \ot \late(x) \bigr) \\
&= \al(y^*) \bigl( \Jh JR(a^*)J \nabh^{1/2}\eta \ot \late(x) \bigr).
\end{align*}
Now $a^* \in \cD(\mu_{-i/2})$ and $R$ and $\mu_t$ commute. So
$R(a^*) \in \cD(\mu_{-i/2})$ and
$\mu_{-i/2}(R(a^*))=R(\mu_{i/2}(a)^*)$. Then we get
$$J R(a^*)J \nabh^{1/2} \subset \nabh^{1/2} J R(\mu_{i/2}(a)^*)
J.$$
Hence we may conclude that $J R(\mu_{i/2}(a)^*)J \eta \in
\cD(\hat{T}^{1/2})$ and
\begin{align*}
(JaJ \ot 1) \Ttil \al(x^*)(\eta \ot \late(y)) &= \al(y^*) \bigl( \hat{T}
J R(\mu_{i/2}(a)^*)J \eta \ot \late(x) \bigr) \\
&= \Ttil \al(x^*) \bigl( J R(\mu_{i/2}(a)^*)J \eta \ot \late(y) \bigr) \\
&= \Ttil (J R(\mu_{i/2}(a)^*)J \ot 1) \al(x^*)(\eta \ot \late(y)).
\end{align*}
Because of the previous lemma we get
\begin{equation}\label{eq*}
(JaJ \ot 1)\Ttil \subset \Ttil (J R(\mu_{i/2}(a)^*)J \ot 1)
\end{equation}
for all $a \in \cD(\mu_{i/2})$. By taking the adjoint we get
$$(J R(\mu_{i/2}(a))J \ot 1) \Ttil^* \subset \Ttil^* (Ja^*J \ot
1)$$
for all $a \in \cD(\mu_{i/2})$. So for all $a \in \cD(\mu_{-i})$
$$(JaJ \ot 1) \nabtil = (JaJ \ot 1) \Ttil^* \Ttil \subset
\Ttil^* (J R (\mu_{-i/2}(a)^*) J \ot 1) \Ttil \subset \nabtil (J
\mu_{-i}(a)J \ot 1).$$
Denoting with $\ga_t$ the automorphism $\Ad \nabtil^{it}$ of $B(H
\ot K)$ we get that for every $a \in \cD(\mu_{-i})$ we have $JaJ
\ot 1 \in \cD(\ga_i)$ and $\ga_i(JaJ \ot 1) = J \mu_{-i}(a)J \ot
1$. Then the results of \cite[4.3 and 4.4]{H4} allow us to
conclude that
$\ga_t(JaJ \ot
1) = J \mu_t(a) J \ot 1$ for every $t \in \R$ and $a \in M$. This
gives
$$(JaJ \ot 1) \nabtil^{-1/2} \subset \nabtil^{-1/2}
(J\mu_{i/2}(a)J \ot 1)$$
for all $a \in \cD(\mu_{i/2})$. Combining this with
equation~\ref{eq*} we get for every $a \in \cD(\mu_{i/2})$
\begin{align*}
(JaJ \ot 1) \Ttil \nabtil^{-1/2} &\subset \Ttil (J
R(\mu_{i/2}(a)^*)J\ot 1) \nabtil^{-1/2} \\
&\subset \Ttil \nabtil^{-1/2}(JR(a^*)J \ot 1) \subset \Jtil(J \Jh
a \Jh J \ot 1).
\end{align*}
So we get
$$(JaJ \ot 1) \Jtil = \Jtil(J \Jh a \Jh J \ot 1) = \Jtil(\Jh J a J
\Jh \ot 1)$$
for every $a \in \cD(\mu_{i/2})$, and hence for every $a \in M$.
Rewriting this we get $(JaJ \ot 1) U = U (J a J \ot 1)$ for every
$a \in M$. This gives $U \in M \ot B(K)$.
\end{proof}
Finally we want to prove that $U$ is a unitary corepresentation of
\qu, namely $(\de \ot \io)(U)=U_{23}U_{13}$. This will be done in
an indirect way in the next section. Nevertheless the results we
use to prove that $U$ is a corepresentation are interesting in
themselves.

\sectie{The unitary implementation is a corepresentation}
The main aim of this section is to prove that the unitary
implementation $U$ is a corepresentation (theorem \ref{44}). On
our way towards the proof of theorem \ref{44} we will solve three
problems which appear naturally in applications (see section
\ref{sectie5} and \cite{K}). First we will see what happens when
we choose a different weight $\te$ on $N$, next we will show how
$U$ changes when the action $\al$ is deformed with an
$\al$-cocycle and finally we will show that in the presence of a
$\sde^{-1}$-invariant weight our implementation agrees with the
one of Enock and Schwartz given by proposition \ref{24}.

In the proof of the first proposition we will make use of
Connes' relative modular theory (see e.g. \cite[3.11,3.12 and
3.16]{S}). When $\te_i$ are \nsf weights on $N$ with GNS-constructions
$(K_i,\pi_i,\la_i)$ $(i=1,2)$, we denote with $J_{2,1}$ the
relative modular conjugation, which is a anti-unitary from $H_1$
to $H_2$. Recall that $J_{1,2}=J_{2,1}^*$. If we denote with $J_i$
the modular conjugation of the weight $\te_i$ we have
$J_{2,1}J_1=J_2 J_{2,1}$ and we denote this unitary with $u$. Then
$u$ is the unique unitary from $K_1$ to $K_2$ which satisfies $u
\pi_1(x) u^* = \pi_2(x)$ for all $x \in N$ and which maps the
positive cone of $K_1$ (determined by the GNS-construction
$(K_1,\pi_1,\la_1)$) onto the positive cone of $K_2$. We will
say that $u$ intertwines the two standard representations of $N$.

Finally we introduce the one-parameter group $\si^{2,1}$ of
isometries of $N$ given by
$$\si^{2,1}_t(x) = [D\te_2 : D \te_1]_t \si^{\te_1}_t(x)$$
for all $x \in N$ and $t \in \R$.
\begin{proposition} \label{41}
Let $\te_i$ be \nsf weights on $N$ with GNS-constructions
$(K_i,\pi_i,\la_i)$ $(i=1,2)$. Let $u$ be the unitary from $K_1$
to $K_2$ intertwining the two standard representations of $N$.
Denote for every $i=1,2$ with $\tetil_i$ the dual weight of
$\te_i$ on $\cros$, with GNS-construction $(H \ot K_i, \io \ot
\pi_i,\latil_i)$. Denote with $U_i \in M \ot B(K_i)$ the unitary
implementation of $\al$ obtained with $\te_i$, as defined in
definition \ref{36}.

Then $1 \ot u$ is the unitary intertwining the two standard
representations of $\cros$. In particular
$$U_2 = (1 \ot u) U_1 (1 \ot u^*).$$
\end{proposition}
\begin{proof}
Let $a \in \Nfih$ and $x \in \cN_{\te_1}$. Let $y \in
\cD(\si^{2,1}_{-i/2})$. Then, by \cite[3.12]{S}, $xy^* \in
\cN_{\te_2}$ and
$$\la_2(xy^*) = J_{2,1} \pi_1(\si^{2,1}_{-i/2}(y)) J_1 \la_1(x).$$
So $(a \ot 1) \al(x) \al(y)^* = (a \ot 1) \al(xy^*) \in
\cN_{\tetil_2}$ and
\begin{align*}
\latil_2 \bigl( (a \ot 1) \al(x) \al(y)^* \bigr) &= \lah(a) \ot \la_2(xy^*)
= \bigl( 1 \ot J_{2,1} \pi_1( \si^{2,1}_{-i/2}(y)) J_1 \bigr) \bigl( \lah(a) \ot
\la_1(x) \bigr) \\  &= \bigl( 1 \ot J_{2,1} \pi_1(\si^{2,1}_{-i/2}(y))
J_1 \bigr) \latil_1 \bigl( (a \ot 1) \al(x) \bigr).
\end{align*}
Because the elements $(a \ot 1) \al(x)$ span a core for $\latil_1$
and because $\latil_2$ is closed (both in the \strong--norm topology), we have for all $z \in \cN_{\tetil_1}$ that $z \al(y)^* \in
\cN_{\tetil_2}$ and
$$\latil_2(z\al(y)^*) = \bigl( 1 \ot J_{2,1} \pi_1(\si^{2,1}_{-i/2}(y))
J_1 \bigr) \latil_1(z).$$
Denoting with $\Jtil_{2,1}$ and $(\sitil_t^{2,1})$ the relative
modular apparatus of the weights $\tetil_2$ and $\tetil_1$, it
follows from \cite[3.12]{S}
that $\al(y) \in \cD(\sitil^{2,1}_{-i/2})$ and
$$\Jtil_{2,1} (\io \ot \pi_1) \bigl( \sitil_{-i/2}^{2,1}(\al(y)) \bigr) \Jtil_1
= 1 \ot J_{2,1} \pi_1(\si^{2,1}_{-i/2}(y)) J_1.$$
Because $[D \tetil_2 : D \tetil_1]_t = \al([D\te_2 :D \te_1]_t)$
for every $t \in \R$ we see that $\sitil^{2,1}_t \na \al = \al
\na \si^{2,1}_t$. So we have
$\sitil_{-i/2}^{2,1}(\al(y))=\al(\si^{2,1}_{-i/2}(y))$. Combining
this with the equation above we get
$$(\io \ot \pi_1)\al \bigl( \si^{2,1}_{-i/2}(y) \bigr) = \Jtil_{2,1}^* (\Jh \ot
J_{2,1}) \bigl( 1 \ot \pi_1(\si^{2,1}_{-i/2}(y)) \bigr) U_1^*.$$
The last formula is valid for all $y \in \cD(\si^{2,1}_{-i/2})$.
Because $U_1$ implements $\al$ we may then conclude that
$U_1=\Jtil_{2,1}^*(\Jh \ot J_{2,1})$.

Then we get
$$1 =U_1 U_1^* = \Jtil_{2,1}^* (\Jh \ot J_{2,1})(\Jh \ot J_1)
\Jtil_1$$
and so $\Jtil_{2,1} \Jtil_1 = 1 \ot J_{2,1} J_1$. Now $u =
J_{2,1}J_1$ and $\Jtil_{2,1} \Jtil_1$ is the unitary intertwining
the two standard representations of $\cros$. This proves the
first claim of the proposition. In particular we get
$$
(1 \ot u)U_1(1 \ot u^*) =(1 \ot u) \Jtil_1 (\Jh \ot J_1) (1 \ot u^*)
= \Jtil_2 (1 \ot u)(\Jh \ot J_1)(1 \ot u^*) = \Jtil_2 (\Jh \ot
J_2) = U_2.
$$
This proves the proposition.
\end{proof}

In the next proposition we will show how the unitary
implementation of an action $\al$ changes when $\al$ is deformed
with an $\al$-cocycle.
\begin{proposition}\label{42}
Let $\al$ be an action of \qu on $N$ and let $\cV \in M \ot N$ be an
$\al$-cocycle in the sense of definition~\ref{15}. Define the
action $\be$ of \qu on $N$ by $\be(x) = \cV \al(x)\cV^*$ for all $x \in
N$. If $\te$ is a \nsf weight on $N$ with GNS-construction $(K,\io,\late)$,
the unitary implementations
$U_\al$ and $U_\be$ of $\al$ and $\be$ obtained with $\te$ satisfy
$$U_\be = \cV \; U_\al \; (\Jh \ot \Jte)\cV^*(\Jh \ot \Jte).$$
In particular $U_\be$ is a corepresentation if and only if $U_\al$
is a corepresentation.
\end{proposition}
\begin{proof}
Because $(\de \ot \io)(\cV) = (1 \ot \cV) (\io \ot \al)(\cV)$ we have
$$(1 \ot \cV^*)(W \ot 1) (1 \ot \cV) = (W \ot 1) (\io \ot \al)(\cV^*).$$
So for every $\xi,\eta \in H$ and with $(e_i)_{i \in I}$ an
orthonormal basis of $H$ we have, by applying $\om_{\xi,\eta} \ot
\io \ot \io$
$$\cV^*(\pih(\om_{\xi,\eta}) \ot 1) \cV = \sum_{i \in I}
(\pih(\om_{e_i,\eta}) \ot 1) \al \bigl( (\om_{\xi,e_i} \ot \io)(\cV^*) \bigr)$$
in the \strong topology. From this it follows that $\cV^*(a \ot 1) \cV
\in \cros$ for all $a \in \Mh$. But $\cV^* \be(x)\cV = \al(x)$ for all
$x \in N$. So
$$\rho: M \kruisje{\be} N \recht \cros : z \mapsto \cV^* z \cV$$
is a well-defined $*$-homomorphism. By symmetry $\rho$ will be
surjective and hence it is a $*$-isomorphism.

Consider now the
dual weights $\tetil_\al$ and $\tetil_\be$ on $\cros$ and $M
\kruisje{\be} N$, with canonical GNS-constructions $(H \ot
K,\io,\latil_\al)$ and $(H \ot K,\io,\latil_\be)$.
Take $\xi \in H$, $b \in \cT_\vfi$ and $x \in \Nte$. Then
$$\cV^*(\pih(\om_{\xi,\la(b)}) \ot 1) \be(x) \cV = \sum_{i \in I}
(\pih(\om_{e_i,\la(b)}) \ot 1) \al \bigl( (\om_{\xi,e_i} \ot \io)(\cV^*)x \bigr)$$
in the \strong topology. For every finite subset $I_0$ of $I$ we
define
$$z_{I_0} := \sum_{i \in I_0} (\pih(\om_{e_i,\la(b)}) \ot 1) \al \bigl( (\om_{\xi,e_i} \ot
\io)(\cV^*)x \bigr).$$
By proposition~\ref{a1} of the appendix we get that $z_{I_0}$ belongs to $\cN_{\tetil_\al}$ and
\begin{align*}
\latil_\al(z_{I_0}) &= \sum_{i \in I_0}
\lah(\pih(\om_{e_i,\la(b)})) \ot (\om_{\xi,e_i} \ot \io)(\cV^*)
\late(x) \\
&= \sum_{i \in I_0} J \si_{i/2}(b) J e_i \ot (\om_{\xi,e_i} \ot
\io)(\cV^*) \late(x) = (J \si_{i/2}(b) J \ot 1)(P_{I_0} \ot 1) \cV^*
(\xi \ot \late(x))
\end{align*}
where $P_{I_0}$ denotes the projection onto $\lspan\{ e_i \mid i
\in I_0\}$. Now define $z:=\cV^*(\pih(\om_{\xi,\la(b)}) \ot 1) \be(x)
\cV$. Then we see that $z_{I_0} \recht z$ \strong and
$$\latil_\al(z_{I_0}) \recht (J \si_{i/2}(b) J \ot 1) \cV^*
(\xi \ot \late(x)) \quad\text{in norm}.$$
So we get that $z \in \cN_{\tetil_\al}$ and
\begin{align*}
\latil_\al(z) &= (J \si_{i/2}(b) J \ot 1) \cV^* (\xi \ot \late(x))
\\ &= \cV^* \bigl( J \si_{i/2}(b)J \xi \ot \late(x) \bigr) = \cV^*
\latil_\be \bigl( (\pih(\om_{\xi,\la(b)}) \ot 1)\be(x) \bigr).
\end{align*}
Because the elements $(\pih(\om_{\xi,\la(b)}) \ot 1) \be(x)$ span
a core for $\latil_\be$ we have $\rho(y) \in \cN_{\tetil_\al}$ for
every $y \in \cN_{\tetil_\be}$ and $\latil_\al(\rho(y)) = \cV^*
\latil_\be(y)$ in that case.

By symmetry $\rho(y) \in \cN_{\tetil_\al}$ if and only if $y \in
\cN_{\tetil_\be}$. But then it is clear that $\Jtil_\be = \cV
\Jtil_\al \cV^*$ and so
$$U_\be = \Jtil_\be(\Jh \ot \Jte) = \cV \Jtil_\al \cV^* (\Jh \ot \Jte)
= \cV U_\al (\Jh \ot \Jte)\cV^*(\Jh \ot \Jte).$$

Now suppose that $U_\al$ is a corepresentation, meaning that $(\de \ot
\io)(U_\al) = U_{\al \, 23} U_{\al \, 13}$. Then
$$(\de \ot \io)(U_\be) = (\de \ot \io)(\cV) \; U_{\al \, 23} U_{\al \, 13} \; (\de \ot
\io)(R \ot L_\te)(\cV)$$
where $L_\te$ is the $*$-anti-isomorphism from $N$ to $N'$ defined by
$L_\te(x)=\Jte x^* \Jte$ for all $x \in N$. Then we can compute
\begin{align*}
(\de \ot \io)(U_\be) &= \cV_{23} (\io \ot \al)(\cV) \; U_{\al \, 23} U_{\al \, 13} \; (R \ot
R \ot L_\te)(\deop \ot \io)(\cV) \\
&= \cV_{23} \; U_{\al \, 23} \cV_{13} U_{\al \, 23}^* \; U_{\al \, 23}  U_{\al \, 13} \; (R \ot R \ot L_\te)(\cV_{13} (\io
\ot \al)(\cV)_{213}) \\
&= \cV_{23} U_{\al \, 23} \cV_{13} U_{\al \, 13} \; (\Jh \ot \Jh \ot \Jte) U_{\al\,
13} \cV_{23}^* U_{\al\, 13}^* \cV_{13}^* (\Jh \ot \Jh \ot \Jte) \\
&= \cV_{23} U_{\al \, 23} \cV_{13} \; (\Jh \ot \Jh \ot \Jte) \cV_{23}^*(\Jh \ot \Jh \ot
\Jte) \; (\Jh \ot \Jh \ot \Jte) U_{\al \, 13}^* \cV_{13}^* (\Jh \ot \Jh \ot \Jte)
\\
&= \cV_{23} U_{\al \, 23} \bigl( (\Jh \ot \Jte)\cV^*(\Jh \ot \Jte) \bigr)_{23} \;\; \cV_{13} U_{\al
\, 13} \bigl( (\Jh \ot \Jte)\cV^*(\Jh \ot \Jte) \bigr)_{13} \\
&= U_{\be \, 23} U_{\be\, 13}.
\end{align*}
So, when $U_\al$ is a corepresentation then $U_\be$ is a corepresentation. By
symmetry also the converse implication holds.
\end{proof}
In proposition~\ref{24} we saw how to construct, with the methods
of Enock and Schwartz, an implementation of an action
$\al$ in the presence of a $\sde^{-1}$-invariant weight. We will show now that
this implementation coincides with the unitary implementation given in
definition~\ref{36}.
\begin{proposition} \label{43}
Let $\al$ be an action of \qu on $N$. Let $\te$ be a \nsf and
$\sde^{-1}$-invariant weight on $N$ with GNS-construction $(K,\io,\late)$. When $V_\te$ is the unitary defined in
proposition~\ref{24} and when $U$ is the unitary implementation of $\al$ defined in
definition~\ref{36}, then $U=V_\te$.
\end{proposition}
\begin{proof}
Recall that
$$(\om_{\xi,\eta} \ot \io)(V_\te) \late(x) = \late \bigl( (\om_{\sde^{1/2}\xi,\eta}
\ot \io) \al(x) \bigr)$$
for all $\xi \in \cD(\sde^{1/2})$ and $\eta \in H$. Because the positive
operators $\sde$ and $\nabh$ strongly commute, we can define the closure $Q$
of the product $\sde \nabh$. Denoting with
$\hogechi_{\mathcal{U}}$ the characteristic function of a subset
$\mathcal{U} \subset \R$ we consider the following subspace of $H$.
$$\cD_0 = \bigcup_{n,m \in \N} \hogechi_{]\frac{1}{n},n[}(\sde)
\hogechi_{]\frac{1}{m},m[}(\nabh) H.$$
Let now $\xi \in H$, $\eta \in \cD_0$, $x \in \cT_\te$ and $y \in \Nte \cap
\Nte^*$. Put again $\hat{T}=\Jh \nabh^{1/2}$ and $\Ttil = \Jtil
\nabtil^{1/2}$. Then
\begin{align*}
(\te_\xi^* \ot 1) \Ttil \al(x^*) (\eta \ot \late(y)) &=(\te_\xi^* \ot
1)\al(y^*)(\hat{T} \eta \ot \late(x)) = \late \bigl( (\om_{\hat{T}\eta,\xi} \ot
\io)\al(y^*)x \bigr) \\
&= \Jte \site_{i/2}(x)^* \Jte \; \late \bigl( (\om_{\hat{T}\eta,\xi} \ot
\io)\al(y^*) \bigr) \\
&= \Jte \site_{i/2}(x)^* \Jte \; (\om_{\sde^{-1/2} \hat{T} \eta,\xi} \ot
\io)(V_\te) \; \late(y^*) \\
&= (\te_\xi^* \ot 1) (1 \ot \Jte \site_{i/2}(x)^* \Jte) \; V_\te
\; (\sde^{-1/2}\hat{T}\eta \ot \late(y^*)).
\end{align*}
Now
$$\sde^{-1/2} \hat{T} \eta = \sde^{-1/2} \Jh \nabh^{1/2} \eta = \Jh \sde^{1/2}
\nabh^{1/2} \eta = \Jh Q^{1/2} \eta.$$
So we may conclude that
$$\Ttil \al(x^*)(\eta \ot \late(y)) = (1 \ot \Jte \site_{i/2}(x)^* \Jte) V_\te
(\Jh \ot \Jte) (Q^{1/2} \ot \nab_\te^{1/2}) (\eta \ot \late(y))$$
for all $\eta \in \cD_0$, $x \in \cT_\te$ and $y \in \Nte \cap \Nte^*$. Because
$\Ttil$ is closed we can conclude that $\eta \ot \late(y) \in \cD(\Ttil)$ and
$$\Ttil (\eta \ot \late(y)) = V_\te
(\Jh \ot \Jte) (Q^{1/2} \ot \nab_\te^{1/2}) (\eta \ot \late(y))$$
for all $\eta \in \cD_0$ and $y \in \Nte \cap \Nte^*$. Because $\cD_0$ is a
core for $Q^{1/2}$ and $\late(\Nte \cap \Nte^*)$ for $\nab_\te^{1/2}$ we get
\begin{equation}\label{inclusie}
V_\te (\Jh \ot \Jte)(Q^{1/2} \ot \nab_\te^{1/2}) \subset \Ttil.
\end{equation}
We now claim that $(Q^{it} \ot \nab_\te^{it}) \Ttil = \Ttil (Q^{it} \ot
\nab_\te^{it})$ for every $t \in \R$. Together with the fact that $\cD(Q^{1/2}
\ot \nab_\te^{1/2}) \subset \cD(\Ttil)$ this leads to the conclusion that $\cD(Q^{1/2}
\ot \nab_\te^{1/2})$ is a core for $\Ttil$. Then we get that the inclusion in
equation~\ref{inclusie} is in fact an equality. Uniqueness of the polar
decomposition gives us $V_\te (\Jh \ot \Jte) = \Jtil$ and so $V_\te = U$.

So we only have to prove our claim. For this choose $x,y \in \Nte$ and $\xi \in
\cD(\hat{T})$. Then using proposition \ref{24} we get
\begin{align*}
(Q^{it} \ot \nab_\te^{it}) \Ttil \al(x^*)(\xi \ot \late(y)) &=(Q^{it} \ot \nab_\te^{it}) \al(y^*)(\hat{T}\xi \ot
\late(x))\\
&= (Q^{it} \ot \nab_\te^{it}) V_\te (1 \ot y^*) V_\te^* (\hat{T} \xi \ot
\late(x))\\
&=V_\te (1 \ot \site_t(y^*)) V_\te^* (Q^{it} \hat{T} \xi \ot \nab_\te^{it}
\late(x))
\end{align*}
because $Q^{it} \ot \nab_\te^{it}$ and $V_\te$ commute. Now observe that
$Q^{it}$ and $\hat{T}$ commute, so that $Q^{it}\xi \in \cD(\hat{T})$ and
\begin{align*}
(Q^{it} \ot \nab_\te^{it}) \Ttil \al(x^*)(\xi \ot \late(y)) &=\al(\site_t(y)^*)
(\hat{T} Q^{it} \xi \ot \late(\site_t(x))) \\
&= \Ttil \al(\site_t(x^*)) (Q^{it} \xi \ot \late(\site_t(y))) \\
&=\Ttil (Q^{it} \ot \nab_\te^{it}) \al(x^*) (\xi \ot \late(y)).
\end{align*}
From this immediately follows our claim, and then the proof of the proposition
is complete.
\end{proof}
With all these results at hand we can now prove the following
theorem.
\begin{theorem} \label{44}
The unitary implementation $U$ of an action $\al$ of \qu on $N$ is
a corepresentation in the sense that $(\de \ot
\io)(U)=U_{23}U_{13}$.
\end{theorem}
\begin{proof}
Consider the bidual action $\alhh$ of $(M,\de)\hoed\Op\hoed\Op$ on
$\Mh \kruisje{\alh} (\cros)$. Let $\te$ be a \nsf weight on $N$
and denote with $\tetiltil$ the bidual weight on $\Mh \kruisje{\alh}
(\cros)$. It follows from proposition~\ref{25} that $\tetiltil$ is a
$J\sde J$-invariant weight for the action $\alhh$. With the
notation of theorem~\ref{26} we define $\rho = \tetiltil \na
\Phi$. Then $\rho$ will be a \nsf and $\sde^{-1}$-invariant weight
on $B(H) \ot N$ for the action $\ga$ of \qu on $B(H) \ot N$.
Combining proposition~\ref{43} and proposition~\ref{24} the
unitary implementation of $\ga$ constructed with the weight $\rho$
is a corepresentation. By proposition~\ref{42} the unitary
implementation of $\mu := (\si \ot \io)(\io \ot \al)$ constructed
with $\rho$ is a corepresentation as well. Then it follows from
proposition~\ref{41} that the unitary implementation $U_\mu$ of
$\mu$ constructed with the \nsf weight $\trace \ot \te$ on $B(H)
\ot N$ will be a corepresentation. Here $\trace$ denotes the usual
trace on $B(H)$.

Represent $N$ on the GNS-space of $\te$ such that $(K,\io,\late)$ is a
GNS-construction for $\te$.
Let $(H_\trace,\pi_\trace,\la_\trace)$ be a GNS-construction for
$\trace$. Then we have a canonical GNS-construction $(H_\trace \ot
K,\pi_\trace \ot \io,\la_{\trace \ot \te})$ for $\trace \ot \te$.
With this we construct the GNS-construction $(H \ot H_\trace \ot
K,\io \ot \pi_\trace \ot \io,\latil_{\trace \ot \te})$ of the dual
weight $(\trace \ot \te)\tildeke$ on $M \kruisje{\mu} (B(H) \ot N)$.
Denote with $\Ttil_{\trace \ot \te}=\Jtil_{\trace \ot \te} \nabtil_{\trace \ot \te}^{1/2}$
the modular operator of this
dual weight. As before we denote with $\Ttil=\Jtil \nabtil^{1/2}$ the modular operator
of the weight $\tetil$ on $\cros$ with GNS-construction $(H \ot K,\io
,\latil)$. It is an easy exercise to check that
$$\Si_{12} \Ttil_{\trace \ot \te} \Si_{12} = J_\trace \ot \Ttil$$
where $\Si_{12}$ flips the first two legs of $H \ot H_\trace \ot
K$. By uniqueness of the polar decomposition we get
$$\Si_{12} \Jtil_{\trace \ot \te} \Si_{12} = J_\trace \ot \Jtil$$
and hence $\Si_{12} U_\mu \Si_{12} = 1 \ot U$. Because $U_\mu$ is
a corepresentation, also $U$ will be a corepresentation.
\end{proof}

\sectie{Subfactors and inclusions of von Neumann algebras} \label{sectie5}
It is well known that there is an important link between
irreducible, depth~2 inclusions of factors and quantum groups.
After a conjecture of Ocneanu the first result in this direction
was proved by David in \cite{Da}, Longo in \cite{Lo} and Szymanski in \cite{Sz}. They
were able to prove that every irreducible, depth~2 inclusion of
$II_1$-factors with finite index has the form $N^\al \subset N$,
where $N$ is a $II_1$-factor and $\al$ is an action of a finite
Kac algebra (i.e. a finite dimensional locally compact quantum
group, or a finite dimensional Hopf $*$-algebra with positive invariant
integral). The restriction on type and index has been removed by
Enock and Nest in \cite{EN} and \cite{E}. There does not appear a
finite quantum group but an arbitrary locally compact quantum
group.

The theory of Enock and Nest is quite technical, but the results
are deep and beautiful. They are important in themselves and serve as a motivation for the concept of a locally compact
quantum group.

Before we describe their result we have to explain a little bit
the basic theory of infinite index inclusions of factors or von
Neumann algebras. So, let us look at an inclusion $N_0 \subset
N_1$ of von Neumann algebras. In this most general setting one can
perform the well known basic construction of Jones. For this we
have to choose a \nsf weight $\te$ on $N_1$ and represent $N_1$ on
the GNS-space of $\te$. Denote with $\Jte$ the modular conjugation
of $\te$. Then we define $N_2=\Jte N_0' \Jte$. Because $N_1' =
\Jte N_1 \Jte$ we have $N_0 \subset N_1 \subset N_2$ and this
inclusion of three von Neumann algebras is called the basic
construction. One can continue in the same way and represent $N_2$
on the GNS-space of some \nsf weight. Then we obtain the von
Neumann algebra $N_3$. Going on we get a tower of von Neumann
algebras
$$N_0 \subset N_1 \subset N_2 \subset N_3 \subset \cdots$$
which is called the Jones tower.

But there is more. In the theory of inclusions of $II_1$-factors
an important role is played by conditional expectations. In the
more general theory being described now, this role will be taken
over by operator valued weights. Before we can explain this, and
also because we need it in the proof of theorem~\ref{52}, we have to explain
Connes' spatial modular theory.
For this we refer to e.g. \cite[\S 7]{S} and \cite[\S III]{Terp}.

Suppose that $N$ is a \vna acting on a Hilbert space $K$. Let $\vfi$ be a
\nsf weight on $N$ and $\psi$ a \nsf weight on $N'$. Let
$(K_\psi,\pi_\psi,\la_\psi)$ be a GNS-construction for $\psi$. For every $\xi
\in K$ we define the densely defined operator $R^\psi(\xi)$ with domain
$\la_\psi(\cN_\psi) \subset K_\psi$ and range in $K$ by $R^\psi(\xi) \la_\psi(x) = x\xi$ for
all $x \in \cN_\psi$. When $\xi \in K$ we can define an operator
$\Te^\psi(\xi)$ in the extended positive part of $B(K)$ by
$$\langle \om_{\eta} , \Te^\psi(\xi) \rangle = \begin{cases} \|
R^\psi(\xi)^* \eta \|^2 &\text{if}\quad \eta \in \cD(R^\psi(\xi)^*) \\
+\infty &\text{else} \end{cases}.$$
In fact $\Te^\psi(\xi) = R^\psi(\xi) R^\psi(\xi)^* + \infty P$ where $P$ is
the projection onto the orthogonal complement of $\cD(R^\psi(\xi)^*)$. Then one
can prove that $\Te^\psi(\xi)$ belongs to $N\Ext$ and it is possible to
define a strictly positive, self-adjoint operator $\displaystyle \frac{d
\vfi}{d \psi}$ on $K$ such that
$$\langle \om_{\xi} , \frac{d \vfi}{d \psi} \rangle = \langle \vfi,
\Te^\psi(\xi) \rangle$$
for all $\xi \in K$. Here we used the extension of the weight $\vfi$ to the
extended positive part of $N$. The operator $\displaystyle \frac{d \vfi}{d
\psi}$ is called the spatial derivative of $\vfi$ with respect to $\psi$.

So, let $N_0 \subset N_1$ be an inclusion of von Neumann algebras
and $T_1$ a \nsf operator valued weight from $N_1$ to $N_0$.
Represent again $N_1$ on the GNS-space of a \nsf weight $\te$.
Let $N_0 \subset N_1 \subset N_2$ be the basic construction. Then
there exists a unique \nsf operator valued weight $T_2$ from
$N_2$ to $N_1$ such that
$$\frac{d (\mu \na T_2)}{d \nu'} = \frac{d \mu}{d \bigl( (\nu \na
T_1)' \bigr)}$$
for all \nsf weights $\mu$ on $N_1$ and $\nu$ on $N_0$. Here we
denote with $\eta'$ the \nsf weight on either $N_2'=\Jte N_0 \Jte$ or
$N_1' = \Jte N_1 \Jte$, given by the formula $\eta'(x)=\eta(\Jte x
\Jte)$ for all positive $x$, whenever $\eta$ is a \nsf weight on
either $N_0$ or $N_1$. The existence of $T_2$ follows from \cite[12.11]{S}.
One can continue in the same way and construct \nsf operator
valued weights $T_i$ from $N_i$ to $N_{i-1}$ anywhere in the Jones
tower.

Next recall that an inclusion of von Neumann algebras $N_0 \subset
N_1$ is said to be
\begin{itemize}
\item irreducible, when $N_1 \cap N_0' = \C$.
\item of depth~2, when $N_1 \cap N_0' \subset N_2 \cap N_0'
\subset N_3 \cap N_0'$ is the basic construction.
\end{itemize}

Finally we describe the notion of regularity as it was introduced
by Enock and Nest in \cite[11.12]{EN}. Let $N_0 \subset N_1$ be an inclusion of
von Neumann algebras. Suppose that $T_1$ is a \nsf operator valued
weight from $N_1$ to $N_0$. Let $N_0 \subset N_1 \subset N_2
\subset N_3 \subset \cdots$ be the Jones tower and construct as
above the operator valued weights $T_2$ from $N_2$ to $N_1$ and
$T_3$ from $N_3$ to $N_2$. Then $T_1$ is called regular when the
restrictions of $T_2$ to $N_2 \cap N_0'$ and of $T_3$ to $N_3 \cap
N_1'$ are both semifinite.

Then we can give the main result of Enock and Nest. Recall that
for a locally compact quantum group \qu we denoted with $(M,\de)'$
the commutant locally compact quantum group, as described in the
introduction.
\begin{theorem}[Enock and Nest] \label{ENest}
Let $N_0 \subset N_1$ be an irreducible, depth~2 inclusion of
factors and let $T_1$ be a regular \nsf operator valued weight
from $N_1$ to $N_0$. Then the von Neumann algebra $M=N_3 \cap
N_1'$ can be given the structure of a locally compact quantum
group $(M,\de)$, such that there exists an outer action $\al$ of
$(M,\de)'$ on $N_1$ satisfying $N_0 = N_1^\al$ and such that the
inclusions $N_0 \subset N_1 \subset N_2$ and $\C \ot N_1^\al \subset
\al(N_1) \subset M' \kruisje{\al} N_1$ are isomorphic.
\end{theorem}
The definition of an outer action is given in definition \ref{outer}.
Further we want to mention that in \cite{E} it is proved
that $(M,\de)$, together with invariant weights and antipode, is
in fact a Woronowicz algebra. But it should be stressed that there
is a small mistake in the proof that the Haar weight is invariant
under the scaling group, so that in fact $(M,\de)$ is an arbitrary
locally compact quantum group, possibly with scaling constant
different from 1.

The main aim of this section is to clarify the conditions of Enock
and Nest's theorem (in particular the regularity condition) and to
prove a converse result: when $\al$ is an integrable and outer
action of \qu on $N$, then the inclusion $N^\al \subset N$ is
irreducible, of depth~2 and the operator valued weight $(\psi \ot
\io)\al$ from $N$ to $N^\al$ is regular.
The same result is stated in \cite[11.14]{EN} for the special case of dual Kac algebra actions, but
the proof is incomplete. The crucial point, our proposition~\ref{54}
identifying two operator valued weights, is not proved in \cite{EN}.
We also remark that it
will follow from corollary~\ref{56} that the actions appearing in Enock and Nest's
theorem are integrable. Further we refer to section~\ref{minimal} for the
link between outer and minimal actions.

First of all we study the following problem. Let $\al$ be an
action of \qu on $N$. Let $N^\al \subset N \subset N_2$ be the
basic construction. When \qu is finite dimensional, it is known
that $N_2$ is a quotient of the crossed product $\cros$ (a proof can be found in \cite[4.1.3]{NSW}, but the result was undoubtedly known before).
More precisely, there exists a surjective $*$-homomorphism $\rho$
from $\cros$ to $N_2$ sending $\al(x)$ to $x$ for all $x \in N$.
So, when $\cros$ is a factor, the inclusion $\C \ot N^\al \subset
\al(N) \subset \cros$ is the basic construction. More
specifically, when $N$ is a $II_1$-factor and $\al$ is an outer
action (or equivalently a free action) of a finite group $G$ on
$N$ it is well known that the crossed product $G \kruisje{\al} N$
can be identified with $(N \cup \{u_t \mid t \in G
\})^{\prime\prime}$, where $N$ is represented standardly and
$(u_t)_{t \in G}$ is the canonical implementation of $\al$. This
can be found in e.g. \cite{JS}.

More generally we look at
the following problem. Suppose that a locally compact group $G$ acts on a
von Neumann algebra $N$ with action $\al$. Then we can construct
the crossed product $G \kruisje{\al} N$ as follows. We represent
$N$ on a Hilbert space $K$ and define operators on $L^2(G) \ot K
\cong L^2(G,K)$ by putting
\begin{alignat*}{2}
(\al(x) \xi)(g) &= \al_{g^{-1}}(x) \xi(g) & & \tekst{for all} g \in G
\tekst{and} x \in N, \xi \in L^2(G,K) \\
(\lambda_g \xi)(h) &= \xi(g^{-1}h) & & \tekst{for all} g,h \in G
\tekst{and} \xi \in L^2(G,K).
\end{alignat*}
Then we define $G \kruisje{\al} N = \bigl( \al(N) \cup \{\lambda_g
\mid g \in G \} \bigr)\dpr$. But, when we represent $N$ standardly
on $K$ and denote with $(u_g)_{g \in G}$ the canonical unitary
implementation of $\al$, we can also define
$$N_2 = \bigl( N \cup \{u_g \mid g \in G \} \bigr)\dpr.$$
Purely algebraicly one would expect to be able to define a
$*$-homomorphism $\rho: G \kruisje{\al} N \recht N_2$ satisfying
$\rho(\al(x)) = x$ for all $x \in N$ and $\rho(\lambda_g)=u_g$ for
all $g \in G$. When the group $G$ is finite, this can be done
easily. In theorem \ref{52} we will prove that the construction of
such a $\rho$ is possible if and only if the action is integrable,
and this will be proved for arbitrary locally compact quantum
group actions. To see the link with the group case, recall that
now the role of the regular representation $(\lambda_g)$ is taken
over by $\lambda(\om) = (\om \ot \io)(W)$ for all $\om \in M_*$.
So we work in fact with the regular representation of the
$L^1$-functions.

Before we come to the proof of our main theorem \ref{52} we characterize the basic
construction $N_2 = \Jte (N^\al)' \Jte$ in terms of the unitary
implementation of $\al$.
\begin{proposition} \label{51}
Let $\al$ be an action of \qu on $N$. Fix a \nsf weight $\te$ on $N$ and let
$N$ act on the GNS-space of $\te$. Let $U$ be the unitary implementation of
$\al$ obtained with $\te$. Let $N_2 = \Jte (N^\al)' \Jte$ be the basic
construction from $N^\al \subset N$. Then
$$N_2 = (N \cup \{(\om \ot \io)(U) \mid \om \in M_* \})^{\prime\prime}.$$
\end{proposition}
\begin{proof}
Because $N_2 = \Jte (N^\al)' \Jte$ we get easily that
$$N_2' = N' \cap \Jte \{ (\om \ot \io)(U^*) \mid \om \in M_* \}' \Jte.$$
But $(\Jh \ot \Jte)U^*(\Jh \ot \Jte) = U$, so that we have
$$N_2' = N' \cap \{(\om \ot \io)(U) \mid \om \in M_* \}'.$$
Because $U$ is a corepresentation the \strong closure of $\{(\om \ot \io)(U) \mid \om \in M_*
\}$ is self-adjoint and then the result follows.
\end{proof}

Then we prove the following result.

\begin{theorem} \label{52}
Let $\al$ be an action of \qu on $N$. Fix a \nsf weight $\te$ on $N$ and let
$N$ act on the GNS-space of $\te$. Let $U$ be the unitary implementation of
$\al$ obtained with $\te$. Let $N_2 = \Jte (N^\al)' \Jte$ be the basic
construction from $N^\al \subset N$. Then the following statements are
equivalent.
\begin{itemize}
\item There exists a normal and surjective $*$-homomorphism $\rho :
\cros \recht N_2$ satisfying
$$\rho(\al(x)) = x \tekst{for all} x \in N \quad\tekst{and}\quad \rho \bigl( (\om \ot
\io)(W) \ot 1 \bigr) = (\om \ot \io)(U^*) \tekst{for all} \om \in M_*.$$
\item The action $\al$ is integrable.
\end{itemize}
\end{theorem}
\begin{proof}[{\bf Proof of the first implication.}]
Let us first suppose the first statement is valid. Because $\Ker \rho$ is a
\strong closed, two-sided ideal of $\cros$ we can take a central projection $P
\in \cros$ such that $\Ker \rho = (\cros)(1-P)$. Let $\rho_P$ be the
restriction of $\rho$ to $(\cros)P$. Then $\rho_P$ is a $*$-isomorphism onto
$N_2$. When $\eta$ is a \nsf weight on $\cros$ we have $\si_t^\eta(P)=P$ for
all $t \in \R$, because $P$ is central. So the restriction $\eta_P$ of $\eta$ to $(\cros)P$ is a \nsf
weight and $\si_t^{\eta_P}$ is the restriction of $\si_t^\eta$ to $(\cros)P$
for all $t \in \R$.

For every \nsf weight $\mu$ on $N$ we define the \nsf weight $\check{\mu}$ on
$N_2$ by $\check{\mu} = (\tilde{\mu})_P \na \rho_P^{-1}$.
Here $\tilde{\mu}$ denotes as before the dual weight on $\cros$.
For every $x \in N$ we have
$$\si_t^{\check{\mu}}(x) =\rho_P \bigl( \si_t^{(\tilde{\mu})_P}(\rho_P^{-1}(x)) \bigr) =
\rho_P \bigl( \si_t^{(\tilde{\mu})_P}(\al(x)P) \bigr) = \rho_P \bigl( \si_t^{\tilde{\mu}}(\al(x))
P \bigr) = \rho_P \bigl( \al(\si_t^\mu(x)) P \bigr) = \si_t^\mu(x).$$
When $\mu$ and $\nu$ are both \nsf weights on $N$ we have
$$[D \check{\mu} : D \check{\nu} ]_t = \rho_P([D (\tilde{\mu})_P : D
(\tilde{\nu})_P ]_t) = \rho_P([D \tilde{\mu} : D \tilde{\nu}]_t P) =
\rho_P(\al([D\mu: D\nu ]_t) P) = [D \mu : D \nu]_t$$
for all $t \in \R$. So, by \cite[12.7]{S}, there exists a unique \nsf operator valued weight
$T_2$ from $N_2$ to $N$ such that $\check{\mu} = \mu \na T_2$ for all \nsf
weights $\mu$ on $N$. So $\mu \na T_2 \na \rho_P = (\tilde{\mu})_P$ for all
\nsf weights $\mu$ on $N$.

When $\nu$ is a \nsf weight on either $N^\al$ or $N$ we denote again with $\nu'$ the
\nsf weight on either $\Jte N^\al \Jte = N_2'$ or $\Jte N \Jte =N'$ given by $\nu'(x)=\nu(\Jte x
\Jte)$ for all positive $x$ in either $N^\al$ or $N$. By \cite[12.11]{S} there exists a
unique \nsf operator valued weight $T_1$ from $N$ to $N^\al$ such that
\begin{equation}\label{operat}
\frac{d(\mu \na T_2)}{d \nu'} = \frac{d \mu}{d \bigl( (\nu \na T_1)' \bigr)}
\end{equation}
for all \nsf weights $\mu$ on $N$ and $\nu$ on $N^\al$.

Choose now a \nsf weight $\te_0$ on $N^\al$. Put $\te_1=\te_0 \na T_1$ and
$\te_2 = \te_1 \na T_2$. When we change the weight $\te$ which was chosen on
$N$ in the beginning of the story, the tower $N^\al \subset N \subset N_2$ will
be transformed into a unitarily equivalent tower. The unitary implementing this
transformation is the unique unitary intertwining the two standard
representations of $N$. This unitary also intertwines the two implementations
of $\al$ by proposition~\ref{41}. Hence also $\rho$ can be transformed.
So we may suppose that $\te=\te_1$.

From equation~\ref{operat} follows that
$$\frac{d\te_2}{d \te_0'} = \frac{d \te_1}{d \te_1'}.$$
So we also have $\dis \frac{d \te_0'}{d\te_2} = \frac{d \te_1'}{d \te_1}$.
But $\dis \frac{d \te_1'}{d \te_1} = \nab_\te^{-1}$ because $K$ is the
GNS-space of $\te=\te_1$. To compute $\dis \frac{d \te_0'}{d\te_2}$ we need a
GNS-construction for the weight $\te_2$. But $\te_2 \na \rho_P = \te \na T_2
\na \rho_P = (\tetil)_P$. So we put $L=P(H \ot K)$ and as before we denote with
$(H \ot K,\io,\latil)$ the GNS-construction of $\tetil$. For every $x \in
\cN_{\te_2}$ we define $\la_{\te_2}(x) = \latil(\rho_P^{-1}(x))$. Then
$\la_{\te_2}(x) \in L$ and it is easy to check that
$(L,\rho_P^{-1},\la_{\te_2})$ is a GNS-construction for $\te_2$. Also observe
that for all $a \in \Nfih$ and $x \in \Nte$ we have $\rho(a \ot 1)x \in
\cN_{\te_2}$ and
$$\la_{\te_2}(\rho(a \ot 1) x) = P(\lah(a) \ot \late(x)).$$
Now choose $z \in \cT_\te$. Then
\begin{equation}\label{starting}
+\infty > \te(z^*z) = \langle \om_{\late(\site_{-i/2}(z))} , \nab_\te^{-1} \rangle
= \langle  \om_{\late(\site_{-i/2}(z))} , \frac{d \te_0'}{d \te_2} \rangle
= \langle \te_0', \Te^{\te_2}\bigl( \late(\site_{-i/2}(z)) \bigr)
\rangle .
\end{equation}
Choose now a family $(\xi_i)_{i \in I}$ of vectors in $K$ such that
$$\te_0'(z) = \sum_{i \in I} \om_{\xi_i}(z)$$
for all $z \in (\Jte N^\al \Jte)^+$. Fix $i \in I$. Then
$$\langle \om_{\xi_i}, \Te^{\te_2} \bigl( \late(\site_{-i/2}(z))
\bigr) \rangle < +\infty$$
and so $\xi_i \in \cD \bigl( R^{\te_2} \bigl( \late(\site_{-i/2}(z)) \bigr)^* \bigr)$. Further
\begin{equation}\label{further}
\langle \om_{\xi_i}, \Te^{\te_2} \bigl(\late(\site_{-i/2}(z))
\bigr) \rangle = \|R^{\te_2} \bigl( \late(\site_{-i/2}(z)) \bigr)^* \xi_i \|^2.
\end{equation}
We will compute the final expression. For this we choose $\om \in \cI$ and $x
\in \Nte$. Recall that the subset $\cI \subset M_*$ was introduced
in the introduction. Observe that
$$R^{\te_2} \bigl( \late(\site_{-i/2}(z)) \bigr)^* \xi_i \in L.$$
So we have
\begin{align*}
\langle R^{\te_2} \bigl( \late(\site_{-i/2}(z)) \bigr)^* \xi_i , \lah \bigl( (\om \ot \io)(W) \bigr) \ot
\late(x) \rangle &= \langle R^{\te_2} \bigl( \late(\site_{-i/2}(z)) \bigr)^* \xi_i, P \bigl( \lah \bigl( (\om \ot \io)(W) \bigr) \ot
\late(x) \bigr) \rangle \\
&= \langle R^{\te_2} \bigl( \late(\site_{-i/2}(z)) \bigr)^* \xi_i, \la_{\te_2} \bigl( \rho \bigl( (\om \ot
\io)(W) \ot 1 \bigr) x \bigr) \rangle \\
&= \langle \xi_i, R^{\te_2} \bigl( \late(\site_{-i/2}(z)) \bigr) \la_{\te_2} \bigl( (\om \ot
\io)(U^*) x \bigr) \rangle \\
&= \langle \xi_i , (\om \ot \io)(U^*)x \late(\site_{-i/2}(z)) \rangle \\
&= \langle \xi_i, (\om \ot \io)(U^*) \Jte z^* \Jte \late(x) \rangle \\
&= \bar{\om} \bigl( (\io \ot \om_{\xi_i,\late(x)}) ((1 \ot \Jte z \Jte) U) \bigr).
\end{align*}
By continuity we get that
$$\langle R^{\te_2} \bigl( \late(\site_{-i/2}(z)) \bigr)^* \xi_i , \lah \bigl( (\om \ot \io)(W) \bigr) \ot
\eta \rangle = \bar{\om} \bigl( (\io \ot \om_{\xi_i,\eta}) ((1 \ot \Jte z \Jte) U)
\bigr)$$
for all $\om \in \cI$, $\eta \in K$ and $z \in \cT_\te$. By proposition~\ref{extra} of the appendix, it follows that
$$(\io \ot \om_{\xi_i,\eta}) \bigl( (1 \ot \Jte z \Jte) U \bigr) \in \Nfi$$
and
$$\la \bigl( (\io \ot \om_{\xi_i,\eta}) ((1 \ot \Jte z \Jte) U) \bigr) = (1 \ot
\te_\eta^*) R^{\te_2} \bigl( \late(\site_{-i/2}(z)) \bigr)^* \xi_i$$
for all $\eta \in K$ and $z \in \cT_\te$. Fix an orthonormal basis $(e_j)_{j
\in J}$ for $K$. Then we may conclude that
\begin{align*}
\|R^{\te_2} \bigl( \late(\site_{-i/2}(z)) \bigr)^* \xi_i \|^2 &=
\sum_{j \in J} \| (1 \ot \te_{e_j}^*) R^{\te_2} \bigl( \late(\site_{-i/2}(z)) \bigr)^*
\xi_i \|^2 \\
&= \sum_{j \in J} \vfi \bigl(
(\io \ot \om_{\xi_i,e_j}) ((1 \ot \Jte z \Jte) U)^*
(\io \ot \om_{\xi_i,e_j}) ((1 \ot \Jte z \Jte) U) \bigr) \\
&= \vfi \bigl( (\io \ot \om_{\xi_i}) (U^*(1 \ot \Jte z^*z \Jte) U) \bigr) \\
&= \vfi( \Jh (\io \ot \om_{\Jte \xi_i}) \al(z^*z) \Jh)
\\
&= \psi \bigl( (\io \ot \om_{\Jte \xi_i}) \al(z^*z) \bigr) \\
&= \langle \om_{\Jte \xi_i} , (\psi \ot \io)\al(z^* z) \rangle.
\end{align*}
Combining this with equation~\ref{further} we get that
$$\langle \om_{\xi_i},
\Te^{\te_2} \bigl( \late(\site_{-i/2}(z)) \bigr) \rangle =
\langle \om_{\Jte \xi_i} , (\psi \ot \io)\al(z^* z) \rangle$$
for all $z \in \cT_\te$ and $i \in I$. Summing over $i$ we get
$$\langle \te_0', \Te^{\te_2} \bigl( \late(\site_{-i/2}(z)) \bigr) \rangle =
\langle \te_0, (\psi \ot \io)\al(z^* z) \rangle$$
for all $z \in \cT_\te$. Using equation~\ref{starting} we get that
$$\te(z^* z) = \langle \te_0, (\psi \ot \io)\al(z^* z) \rangle$$
for all $z \in \cT_\te$. Hence the normal faithful weight $\te_0 \na (\psi \ot
\io)\al$ is semifinite. From \cite[11.7]{S} it follows that $(\psi \ot
\io) \al$ is semifinite. So $\al$ is integrable.

{\bf Proof of the second implication.} The second implication can be proved along
the same lines as in the case of an Abelian group action, see \cite{P}.
So let us suppose that $\al$ is
integrable. Choose a \nsf weight $\te_0$ on $N^\al$ and put $\te = \te_0 \na
(\psi \ot \io)\al$. Then $\te$ is a \nsf weight on $N$. Represent $N$ on the
GNS-space $K$ of $\te$ such that $(K,\io,\late)$ is a GNS-construction for
$\te$. Choose a family of vectors $(\xi_i)_{i \in I}$ in $K$ such that
$$\te_0(x) = \sum_{i \in I} \om_{\xi_i}(x) \tekst{for all} x \in
(N^\al)^+.$$
Define $L = \bigoplus_{i \in I} H \ot K$ and let $\pi$ be the $I$-fold direct
sum of the identical representation $\io$ of $\cros$ on $H \ot K$.
Recall that for any operator valued weight $T$ we define $\cN_T$
as the left ideal of elements $x$ for which $T(x^* x)$ is bounded.
Also recall that we introduced the canonical GNS-construction
$(H,\io,\Gamma)$ for $\psi$ in the introduction. When $z \in
\cN_{\psi \ot \io}$ we define $(\Ga \ot \io)(z) \in B(K,H \ot K)$ by $(\Ga \ot
\io)(z) \late(x) = (\Ga \ot \late)(z(1 \ot x))$ for all $x \in \Nte$,
where $\Ga \ot \late$ denotes the canonical GNS-map of $\psi \ot \te$.
One can check easily that $(\Ga \ot \io)(z)^* (\Ga \ot \io)(z) =
(\psi \ot \io)(z^*z)$. For this see e.g. \cite[10.6]{EN}. Put $T =
(\psi \ot \io) \al$. For all $x \in \cN_T \cap \Nte$ we define
$$\cV \late(x) = \bigoplus_{i \in I} (\Ga \ot \io)\al(x) \xi_i.$$
Observe that $\cV$ is well-defined:
$$\sum_{i \in I} \| (\Ga \ot \io)\al(x) \xi_i \|^2 = \sum_{i \in I}
\om_{\xi_i} \bigl( (\psi \ot \io)\al(x^* x) \bigr) = \langle \te_0,(\psi \ot
\io)\al(x^* x) \rangle = \te(x^* x) < \infty.$$
Because $\cN_T \cap \Nte$ is a \strong--norm core for $\late$ we get that
$\late(\cN_T \cap \Nte)$ is dense in $K$. So $\cV$ can be extended uniquely to an
isometry from $K$ to $L$.

We now want to prove that the range of $\cV$ is invariant under $\pi(\cros)$.
So we first choose $y \in N$. Then for every $x \in \cN_T \cap \Nte$ we have
\begin{align*}
\pi(\al(y)) \cV \late(x) &= \bigoplus_{i \in I} \al(y) (\Ga \ot \io)\al(x) \xi_i
\\ &= \bigoplus_{i \in I} (\Ga \ot \io)\al(yx) \xi_i = \cV \late(yx) = \cV y
\late(x).
\end{align*}
Next we will look at the invariance under $\pi(\Mh \ot \C)$.
Analogously as in proposition~\ref{a2} of the appendix we have that for every $x \in \Npsi$,
$\xi \in \cD(\sde^{1/2})$ and $\eta \in H$, $(\om_{\sde^{1/2} \xi,\eta} \ot
\io)\de(x) \in \Npsi$ and
$$\Ga \bigl( (\om_{\sde^{1/2} \xi,\eta} \ot \io)\de(x) \bigr) = (\om_{\xi,\eta} \ot
\io)(W^*) \Ga(x).$$
Then it follows easily that for all $z \in \cN_{\psi \ot \io}$ we have
$$x:= (\om_{\sde^{1/2} \xi,\eta} \ot \io \ot \io)(\de \ot \io)(z) \in \cN_{\psi
\ot \io} \tekst{and} (\Ga \ot \io)(x) = \bigl( (\om_{\xi,\eta} \ot \io)(W^*) \ot 1
\bigr) (\Ga \ot \io)(z).$$
So for all $\xi \in \cD(\sde^{1/2})$, $\eta \in H$ and $x \in \cN_T \cap \Nte$ we have
\begin{align*}
\pi \bigl( (\om_{\xi,\eta} \ot \io)(W^*) \ot 1 \bigr) \cV \late(x) &= \bigoplus_{i \in I}
\bigl( (\om_{\xi,\eta} \ot \io)(W^*) \ot 1 \bigr)(\Ga \ot \io) \al(x) \xi_i \\
&= \bigoplus_{i \in I} (\Ga \ot \io) \bigl( (\om_{\sde^{1/2} \xi,\eta} \ot \io \ot \io)(\de \ot
\io)\al(x) \bigr) \xi_i \\
&= \bigoplus_{i \in I} (\Ga \ot \io) \bigl( \al( (\om_{\sde^{1/2} \xi,\eta} \ot
\io)\al(x) ) \bigr) \xi_i \\
&= \cV \late \bigl( (\om_{\sde^{1/2} \xi,\eta} \ot
\io)\al(x) \bigr) = \cV(\om_{\xi,\eta} \ot \io)(U) \late(x)
\end{align*}
by propositions~\ref{24} and \ref{43}. So the range of $\cV$ is invariant under
$\pi(\cros)$. Then we can define a $*$-homomorphism
$$\rho: \cros \recht B(K) : \rho(z) = \cV^* \pi(z) \cV.$$
By the computations above we get that $\rho(\al(x)) = x$ for all $x \in N$ and
$(\io \ot \rho)(W \ot 1) = U^*$. Then it follows from proposition~\ref{51} that
$\rho(\cros) = N_2$ and so the theorem is proved.
\end{proof}

One can also prove the following more general kind of result,
where we do not specify what $\rho$ should be.
\begin{proposition} \label{53}
Let $\al$ be an action of \qu on $N$. Fix a \nsf weight $\te$ on $N$ and let
$N$ act on the GNS-space $K$ of $\te$. Consider the inclusions
$$\C \ot N^\al \subset \al(N) \subset \cros \tekst{and} N^\al \subset N
\subset N_2=\Jte (N^\al)' \Jte.$$
Then the following statements are equivalent.
\begin{itemize}
\item There exists a surjective $*$-homomorphism $\rho$ from $\cros$ to $N_2$ such
that $\rho$ is an isomorphism of $\al(N)$ onto $N$ and of $\C \ot
N^\al$ onto $N^\al$.
\item The action $\al$ is cocycle-equivalent with an integrable action $\be$
satisfying $N^\be=N^\al$.
\end{itemize}
\end{proposition}
\begin{proof}[{\bf Proof of the first implication.}]
Suppose the first statement is true. Because $N$ is represented
on the GNS-space of $\te$,
there exists a unitary $u$ on $K$ such that $\rho(\al(x)) = u x u^*$ for all $x
\in N$ and $u \Jte = \Jte u$. Define $\tilde{\rho}$ from $\cros$ to $B(K)$ by
$\tilde{\rho}(z) = u^* \rho(z)u$ for all $z \in \cros$. Then
$\tilde{\rho}(\al(x)) = x$ for all $x \in N$. Further
$$u^* N^\al u = u^* \rho(\C \ot N^\al)u = u^*\rho(\al(N^\al))u= N^\al.$$
Because $u \Jte = \Jte u$, we get $u^* N_2 u=N_2$, which leads to $\tilde{\rho}(\cros)=N_2$.

So we may suppose from the beginning that $\rho(\al(x))=x$ for all $x \in N$.
Define the unitary $X \in M \ot B(K)$ by
$$X = (\Jh \ot \Jte) (\io \ot \rho)(W \ot 1)(\Jh \ot \Jte).$$
Put $\cV = XU^*$. Clearly $\cV \in M \ot B(K)$ and
$$(\Jh \ot \Jte) \cV (\Jh \ot \Jte) = (\io \ot \rho)(W \ot 1) U.$$
For every $x \in N$ we have
\begin{align*}
(\io \ot \rho)(W \ot 1) U (1 \ot x) &= (\io \ot \rho)(W \ot 1) \al(x) U = (\io
\ot \rho) \bigl( (W \ot 1)(\io \ot \al)\al(x) \bigr) U \\
&= (\io \ot \rho) \bigl( (1 \ot \al(x))(W \ot 1) \bigr) U = (1 \ot x) (\io \ot \rho)(W \ot
1)U.
\end{align*}
So we get $(\io \ot \rho)(W \ot 1) U \in M \ot N'$ and hence $\cV \in M \ot N$.
In the next computation we denote again with
$L_\te$ the $*$-anti-automorphism of $B(K)$ given by $L_\te(x)=\Jte x^* \Jte$
for all $x \in B(K)$. Then we have
\begin{align*}
(\de \ot \io)(\cV) &= (\de \ot \io)(R \ot L_\te)(\io \ot \rho)(W^* \ot 1) \;
(\de \ot \io)(U^*) = (R \ot R \ot L_\te)(\deop \ot \rho)(W^* \ot 1) \; U^*_{13}
U^*_{23} \\
&= (R \ot R \ot L_\te)(\io \ot \io \ot \rho)(W^*_{13} W^*_{23}) \; U^*_{13}
U^*_{23} \\ &=  \bigl( (R \ot L_\te)(\io \ot \rho)(W^* \ot 1) \bigr)_{23} U_{23}^* \; U_{23} \bigl( (R
\ot L_\te)(\io \ot \rho)(W^* \ot 1) U^* \bigr)_{13} U_{23}^* \\
&= \cV_{23} (\io \ot \al)(\cV).
\end{align*}
So $\cV$ is a $\al$-cocycle. Define the action $\be$ of \qu on $N$ given by
$\be(x) = \cV \al(x) \cV^*$ for all $x \in N$. Then $\be(x) = X (1 \ot x) X^*$
for all $x \in N$. Because the \strong closure of $\{(\om \ot \io)(X) \mid \om
\in M_*\}$ equals $\Jte \rho(\Mh \ot \C) \Jte$ we get that $N^\be = \Jte
\rho(\cros)'\Jte = N^\al$.

To conclude the proof of the first implication we have to show that $\be$ is
integrable. For this we will use the previous theorem. From
proposition~\ref{42} it follows that the unitary implementation $U_\be$ of
$\be$ is given by
$$U_\be = \cV \; U \; (\Jh \ot \Jte)\cV^*(\Jh \ot \Jte)= \cV (\io \ot \rho)(W^* \ot
1).$$
From the proof of proposition~\ref{42} we also get that $z \mapsto \cV^* z \cV$
gives an isomorphism from $M \kruisje{\be} N$ onto $\cros$. So we can define
$$\tilde{\rho} : M \kruisje{\be} N \recht N_2 : \tilde{\rho}(z)=\rho(\cV^* z
\cV).$$
Then $\tilde{\rho}$ is a surjective $*$-homomorphism onto $N_2 = \Jte (N^\al)'
\Jte = \Jte (N^\be)' \Jte$ and clearly $\tilde{\rho}(\beta(x))=x$ for all $x \in N$.
Because $\cV$ is an $\al$-cocycle we get that
$$(1 \ot \cV^*)(W^* \ot 1)(1 \ot \cV) = (\io \ot \al)(\cV)(W^* \ot 1).$$
From this it follows that
$$(\io \ot \tilde{\rho})(W^* \ot 1) = (\io \ot \rho) \bigl( (1 \ot \cV^*)(W^* \ot 1)(1 \ot
\cV) \bigr) = \cV (\io \ot \rho)(W^* \ot 1) = U_\be.$$
By the previous theorem we get that $\be$ is integrable.

{\bf Proof of the second implication.} Conversely suppose that the second
statement is valid and take such an action $\be$. Let $\cV$ be an $\al$-cocycle
such that $\be(x)=\cV \al(x) \cV^*$ for all $x \in N$. It follows from the
proof of proposition~\ref{42} that
$$\Phi: \cros \recht M \kruisje{\be} N : \Phi(z) = \cV z \cV^*$$
is an isomorphism and $\Phi(\al(x)) = \be(x)$ for all $x \in N$. By the
previous theorem we can find a surjective $*$-homomorphism $\tilde{\rho}$ from
$M \kruisje{\be} N$ onto $\Jte (N^\be)' \Jte$ satisfying $\tilde{\rho}(\be(x))
= x$ for all $x \in N$. Putting $\rho=\tilde{\rho} \na \Phi$ and observing that
$N^\al=N^\be$ we get the first statement.
\end{proof}
We do not know an example of a non-integrable action $\al$ which
is cocycle-equivalent with an integrable action $\be$ satisfying
$N^\al=N^\be$, but it seems to be natural that this kind of
actions will exist. We will now specify a case in which it cannot
exist. This should be compared with the example of a finite group
acting outerly on a factor as described above.
\begin{definition} \label{outer}
An action $\al$ of a locally compact quantum group \qu on $N$ is
called outer when
$$\cros \cap \al(N)' = \C.$$
\end{definition}
\begin{corollary} \label{56}
Let $\al$ be an outer action of \qu on $N$. Choose again a \nsf
weight $\te$ on $N$ and represent $N$ on the GNS-space of $\te$.
Let $N_2 = \Jte (N^\al)' \Jte$ be the basic construction. Then the
inclusions
$$\C \ot N^\al \subset \al(N) \subset \cros \tekst{and} N^\al \subset N
\subset N_2$$
are isomorphic if and only if $\al$ is integrable.
\end{corollary}
\begin{proof}
When $\al$ is integrable, one can use theorem~\ref{52} and then
observe that the $*$-homomorphism $\rho$ is faithful because
$\cros$ is a factor.

Next suppose that the inclusions stated above are isomorphic. By
proposition~\ref{53} there exists an integrable action $\be$ which
is cocycle equivalent with $\al$ and satisfies $N^\be=N^\al$. Let
$\cV \in M \ot N$ be an $\al$-cocycle such that $\be(x)=\cV\al(x) \cV^*$
for all $x \in N$. Then for all $x \in N^\al= N^\be$ we have
$$1 \ot x = \be(x) = \cV\al(x) \cV^* = \cV (1 \ot x) \cV^*.$$
Hence we get $\cV \in M \ot (N \cap (N^\al)')$. From our assumption
and the fact that $\al$ is outer it follows that $N_2 \cap N'=\C$.
But then also
$$(N^\al)' \cap N = \Jte (N_2 \cap N') \Jte = \C.$$
So we can take $u \in M$ such that $\cV = u \ot 1$. Because $\cV$ is
an $\al$-cocycle we get that $\de(u)=u \ot u$. By the unicity of right
invariant weights on \qu there exists a number $\lambda > 0$ such
that $\psi(u^* a u) = \sla \psi(a)$ for all $a \in M^+$. Then we
get that for all $x \in N^+$ we have $(\psi \ot \io) \al(x)= \sla
(\psi \ot \io) \be(x)$. Because $\be$ is integrable it follows
that $\al$ is integrable.
\end{proof}
There exist outer actions which are not integrable: see \ref{example}.
Combining the previous result with theorem \ref{ENest} we get that all
actions coming out of Enock and Nest's construction are
integrable.

Next we turn towards the notion of a regular operator valued
weight. Suppose $\al$ is an integrable action of \qu on $N$ and
suppose that the $*$-homomorphism $\rho$ given by theorem~\ref{52} is
faithful. This will of course be the case whenever $\cros$ is a factor,
but also when $\al$ is a dual action or a semidual action. The latter follows from proposition~\ref{semidual}.
Then we can prove that the operator valued weight $(\psi \ot
\io)\al$ from $N$ to $N^\al$ is regular. More precisely, we will
do the following. By our assumption the basic construction $N^\al
\subset N \subset N_2$ is isomorphic with $\C \ot N^\al \subset
\al(N) \subset \cros$ through the isomorphism $\rho$.
Let us denote with $T_1$ the operator valued weight $(\psi \ot \io)\al$ from $N$ to
$N^\al$. Then we can
construct the operator valued weight $T_2$ from $N_2$ to $N$ by modular theory and the basic construction, as
described above.
Through the isomorphism $\rho$ the operator valued weight
$T_2$ is transformed to an operator valued weight from $\cros$ to $\al(N)$. In
the next proposition we prove that this operator valued weight is equal to the
canonical operator valued weight $T = (\vfih \ot \io \ot \io) \alh$ from
$\cros$ to $\al(N)$.
\begin{proposition} \label{54}
Let $\al$ be an integrable action of \qu on $N$. Suppose that the
$*$-homomorphism $\rho$ constructed in theorem~\ref{52} is faithful. Denote
with $T_2$ and $T$ the operator valued weights defined above. Then $\rho \na T
= T_2 \na \rho$.
\end{proposition}
\begin{proof}
For clarity we stress that $T_1$ is the operator valued weight
$(\psi \ot \io)\al$ from $N$ to $N^\al$, that $T_2$ is obtained
out of $T_1$ by modular theory and the basic construction, and it
goes from $N_2$ to $N$. Finally $T$ is the canonical operator
valued weight $(\vfih \ot \io \ot \io) \alh$ from $\cros$ to
$\al(N)$, giving the dual weights by the formula $\tetil = \te \na
\al^{-1} \na T$ for all \nsf weights $\te$ on $N$.

Choose a \nsf weight $\te_0$ on $N^\al$. Put $\te = \te_0 \na T_1$ and let
$\te_2 = \tilde{\te} \na \rho^{-1}$.
We will prove that $\te_2= \te \na T_2$.
As in the proof of theorem~\ref{52} we may
suppose that $N$ is represented on the GNS-space of $\te$ such that
$(K,\io,\late)$ is a GNS-construction for $\te$. Let $(H \ot K,\io,\latil)$ be
the canonical GNS-construction for $\tetil$  and put $\la_{\te_2} = \latil \na
\rho^{-1}$. We now make a kind of converse reasoning of the proof
of theorem~\ref{52}.
Denote again with $\te'_0$ the \nsf weight on $\Jte N^\al \Jte = N'_2$ given by
$\te'_0(x) = \te_0(\Jte x \Jte)$ for all positive $x$.
We claim that for all $z \in \cT_\te$
\begin{equation} \label{claim}
\langle \om_{\late(\site_{-i/2}(z))} , \frac{d \te_0'}{d \te_2} \rangle = \langle \te_0 ,
(\psi \ot \io) \al(z^* z) \rangle.
\end{equation}
So choose $z \in \cT_\te$. Take a family of vectors $(\xi_i)_{i
\in I}$ in $K$ such that
$$\te_0(x) = \sum_{i \in I} \langle x \Jte \xi_i,\Jte \xi_i
\rangle$$
for all $x \in (N^\al)^+$. Because $\langle \te_0, (\psi \ot
\io)\al(z^* z) \rangle = \te(z^* z) < \infty$ we have
$$\langle \om_{\Jte \xi_i}, (\psi \ot
\io)\al(z^* z) \rangle < \infty$$
for all $i \in I$. Fix $i \in I$. Then we conclude from the
previous formula that
$$\vfi \bigl( (\io \ot \om_{\xi_i})(U^* (1 \ot \Jte z^* z \Jte) U) \bigr) <
\infty.$$
So, when $(e_j)_{j \in J}$ is an orthonormal basis for $K$ we can
define the element $\eta \in H \ot K$ by
$$\eta := \sum_{j \in J} \la \bigl( (\io \ot \om_{\xi_i,e_j})( (1 \ot
\Jte z \Jte) U) \bigr) \ot e_j.$$
It is easy to check that for all $\mu \in K$ we have $(\io \ot
\om_{\xi_i,\mu}) ( (1 \ot \Jte z \Jte) U ) \in \Nfi$ and
$$\la \bigl( (\io \ot
\om_{\xi_i,\mu}) ( (1 \ot \Jte z \Jte) U ) \bigr) = (1 \ot
\te_\mu^*) \eta.$$
Using the notation $\cI \subset M_*$ introduced in the
introduction, we get for all $\om \in \cI$ and $x \in \Nte$ that
\begin{align*}
\langle \xi_i, R^{\te_2} \bigl( \late(\site_{-i/2}(z)) \bigr) \la_{\te_2} \bigl(
(\om \ot \io)(U^*) x \bigr) \rangle
&= \langle \xi_i , (\om \ot \io)(U^*) \Jte z^* \Jte \late(x)
\rangle \\
&= \bar{\om} \bigl( (\io \ot \om_{\xi_i,\late(x)}) ( (1 \ot \Jte z
\Jte) U) \bigr) \\
&= \langle \eta, \lah \bigl( (\om \ot \io) (W) \bigr) \ot \late(x)
\rangle \\
&= \langle \eta, \la_{\te_2} \bigl( (\om \ot \io)(U^*)x \bigr) \rangle.
\end{align*}
From this we get that
$$\langle \xi_i, R^{\te_2} \bigl( \late(\site_{-i/2}(z)) \bigr)
\la_{\te_2}(y) \rangle = \langle \eta,\la_{\te_2}(y) \rangle$$
for all $y \in \cN_{\te_2}$. Hence $\xi_i \in \cD \bigl( R^{\te_2} \bigl( \late(\site_{-i/2}(z))
\bigr)^* \bigr)$ and
$$\| R^{\te_2} \bigl( \late(\site_{-i/2}(z)) \bigr)^* \xi_i \|^2 = \|\eta\|^2
= \langle \om_{\Jte \xi_i}, (\psi \ot \io)\al(z^* z)
\rangle.$$
This means that
$$\langle \om_{\xi_i},
\Theta^{\te_2} \bigl( \late(\site_{-i/2}(z)) \bigr)
\rangle = \langle \om_{\Jte \xi_i}, (\psi \ot \io)\al(z^* z)
\rangle.$$
Summing over $i$ we get our claim stated in equation~\ref{claim}.
But now
$$\langle \te_0 ,
(\psi \ot \io) \al(z^* z) \rangle = \te(z^* z) = \langle
\om_{\late(\site_{-i/2}(z))} , \nab_\te^{-1} \rangle,$$
and so
$$\langle \om_{\late(\site_{-i/2}(z))} , \nab_\te^{-1} \rangle =
\langle \om_{\late(\site_{-i/2}(z))} , \frac{d \te'_0}{d \te_2}
\rangle$$
for all $z \in \cT_\te$.
Next we claim that $\dis \frac{d \te_0'}{d \te_2}$ and $\nab_\te^{-1}$ commute
strongly. Then we will be able to conclude that $\dis \frac{d \te_0'}{d
\te_2}=\nab_\te^{-1}$. But then $\dis \frac{d \te_2}{d \te_0'} =
\nab_\te$, and so we will get
$$\frac{d(\te \na T_2)}{d \te_0'} = \frac{d \te}{d \bigl( (\te_0 \na T_1)' \bigr)} = \frac{d
\te}{d \te'} = \nab_\te= \frac{d \te_2}{d \te_0'}.$$
So we may conclude that $\te_2=\te \na T_2$.
By definition of $\tetil$ we have $\tetil = \te \na \rho \na T$ and then $\te \na \rho \na T =
\tetil = \te_2 \na \rho = \te \na T_2 \na \rho$. By \cite[11.13]{S} we get that $\rho \na T
= T_2 \na \rho$.

So we only have to prove our claim. Hence we want to prove that $\dis \frac{d \te_0'}{d
\te_2}$ and $\nab_\te^{it}$ commute for every $t \in \R$. For this it is
sufficient to prove that $\Ad \nab_\te^{it}$ leaves both $\Jte N^\al \Jte$ and
$N_2$ invariant and
$$\te_0' \na \Ad \nab_\te^{it} = \te_0' \quad\tekst{and}\quad \te_2 \na \Ad
\nab_\te^{it} = \te_2$$
for all $t \in \R$. When $x \in N^\al$ we have
$$\nab_\te^{it} \Jte x \Jte \nab_\te^{-it} = \Jte \site_t(x) \Jte = \Jte
\si^{\te_0}_t(x) \Jte \in \Jte N^\al \Jte.$$
Then it is immediately clear that $\te_0' \na \Ad \nab_\te^{it} = \te_0'$.

Because $N_2 = (\Jte N^\al \Jte)'$ we have that $\Ad \nab_\te^{it}$ leaves
$N_2$ invariant. Recall that we denoted with $(\sitil_t)$ the
modular group of $\tetil$ on $\cros$.
Then we have, for all $x \in N$
\begin{equation}\label{second}
\nab_\te^{it} \rho(\al(x)) \nab_\te^{-it} = \nab_\te^{it} x \nab_\te^{-it}
=\site_t(x) =\rho \bigl( \al(\site_t(x)) \bigr)= \rho \bigl( \sitil_t(\al(x)) \bigr).
\end{equation}
Finally, for all $\om \in B(H)_*$ we have by proposition~\ref{43} and \ref{24}
that
\begin{align*}
\nab_\te^{it} \rho \bigl( (\om \ot \io)(W) \ot 1 \bigr) \nab_\te^{-it} &= \nab_\te^{it} (\om
\ot \io)(U^*) \nab_\te^{-it} = (Q^{it} \om Q^{-it} \ot \io)(U^*) \\
&= \rho \bigl( (Q^{it} \om Q^{-it} \ot \io)(W) \ot 1 \bigr) \\ &= \rho \bigl( (Q^{it} \ot
\nab_\te^{it}) \bigl( (\om \ot \io)(W) \ot 1 \bigr) (Q^{-it} \ot \nab_\te^{-it})
\bigr)
\end{align*}
where we used that $W(Q^{it} \ot Q^{it})=(Q^{it} \ot Q^{it})W$. From the proof
of proposition~\ref{43} it follows that $\nabtil^{it} = Q^{it} \ot
\nab_\te^{it}$ and so we see that
$$\nab_\te^{it} \rho(a \ot 1) \nab_\te^{-it} = \rho(\sitil_t(a \ot 1))$$
for all $a \in \Mh$ and $t \in R$. Combining this with equation~\ref{second} we
get that $\nab_\te^{it} \rho(z) \nab_\te^{-it} = \rho(\sitil_t(z))$ for all $z
\in \cros$ and $t \in \R$. Then we get immediately that $\te_2 \na \Ad
\nab_\te^{it} = \te_2$ for all $t \in \R$.

This proves our claim and ends the proof of the proposition.
\end{proof}
\begin{corollary} \label{regular}
Under the same assumptions as in proposition~\ref{54}, the
operator valued weight $(\psi \ot \io) \al$ from $N$ to $N^\al$ is
regular.
\end{corollary}
\begin{proof}
Using the notations introduced above we will identify the
inclusions $N^\al \subset N \subset N_2$ and $\C \ot N^\al \subset
\al(N) \subset \cros$. Then we get that $T_2 = (\vfih \ot \io \ot
\io) \alh$. Now it is obvious that $\Mh \ot \C \subset \cros \cap
(\C \ot N^\al)'$ and $\cN_{\vfih} \ot \C \subset \cN_{T_2}$. So
the restriction of $T_2$ to $N_2 \cap (N^\al)'$ is semifinite.

Next observe that $\al(N)= (\cros)^{\alh}$. Applying the first
part of the proof to the dual action $\alh$,
which is integrable and for which the $*$-homomorphism $\rho$ is
faithful by proposition \ref{semidual}, we get that the
restriction of $T_3$ to $N_3 \cap N_1'$ is semifinite.
\end{proof}

As a final ingredient for the converse of Enock and Nest's theorem we look at
depth~2 inclusions. The assumption of the following proposition may seem
strange, but one can immediately look at the corollary for a more clear result.
\begin{proposition} \label{depth}
Let $\al$ be an action of \qu on $N$ such that $\C \ot N^\al \subset \al(N) \subset
\cros$ is the basic construction. Then the inclusion $N^\al \subset N$ has
depth~2.
\end{proposition}
\begin{proof}
Choose a \nsf weight $\te$ on $N$ and let $\tetil$ be the dual weight on
$\cros$. Represent $N$ on the GNS-space of $\te$ such that $(K,\io,\late)$ is a
GNS-construction for $\te$. Let $(H \ot K,\io,\latil)$ be the canonical
GNS-construction for $\tetil$ and denote with $\Jtil$ the modular conjugation
of $\tetil$. Then it follows from definition~\ref{36} that $U=\Jtil (\Jh \ot \Jte)$
is the unitary implementation of $\al$. The basic construction from $\al(N)
\subset \cros$ is then given by
$$\Jtil \al(N)' \Jtil = \Jtil U (B(H) \ot N') U^* \Jtil = B(H) \ot N.$$
To prove that $N^\al \subset N$ has depth~2, we have to show that
$$\al(N \cap (N^\al)') \subset (\cros) \cap (\C \ot N^\al)' \subset B(H) \ot (N
\cap (N^\al)')$$
is the basic construction. But it is immediately clear that the restriction of
$\al$ to $N \cap (N^\al)'$ is an action $\be$ of \qu on $N \cap (N^\al)'$. So
by the first part of the proof it is sufficient to prove that
\begin{equation}\label{kruis}
M \kruisje{\be} (N \cap (N^\al)') = (\cros) \cap (\C \ot N^\al)'.
\end{equation}
Now it follows from theorems~\ref{26} and \ref{27} that
\begin{align*}
M \kruisje{\be} (N \cap (N^\al)') &= \{ z \in B(H) \ot (N \cap (N^\al)') \mid
(\io \ot \be)(z) = V_{12} z_{13} V_{12}^* \} \\
\intertext{and}
\cros &= \{z \in B(H) \ot N \mid (\io \ot \al)(z)  = V_{12} z_{13} V_{12}^* \}.
\end{align*}
From this we can immediately deduce equation~\ref{kruis}, and that concludes
the proof.
\end{proof}
Although the following result is an immediate corollary of the
previous one, we include it for completeness. The first statement
is clear and the next two statements follow from the first, using
proposition \ref{semidual} for the last one.
\begin{corollary} \label{cor510}
Let $\al$ be an action of \qu on $N$.
\begin{itemize}
\item If $\al$ is integrable and the $*$-homomorphism in theorem \ref{52} is
faithful, then the inclusion $N^\al \subset N$ has depth~2.
\item If $\al$ is integrable and $\cros$ is a factor, then the inclusion $N^\al
\subset N$ has depth~2.
\item The inclusion $\al(N) \subset \cros$ has depth~2.
\end{itemize}
\end{corollary}
We now prove the announced result giving a converse to the theorem of Enock and
Nest.
\begin{proposition}
Let $\al$ be an integrable outer action of \qu on $N$. Then the operator valued
weight $(\psi \ot \io)\al$ from $N$ to $N^\al$ is regular. Further the
inclusion $N^\al \subset N$ is irreducible and has depth~2.
\end{proposition}
\begin{proof}
Because $\cros$ is a factor the $*$-homomorphism $\rho$ from theorem \ref{52} is
faithful. Then we apply corollary~\ref{regular} to obtain the regularity of $(\psi \ot \io)\al$ and
corollary \ref{cor510} to get that $N^\al \subset N$ has depth~2. It is clear that $N^\al \subset
N$ is irreducible, because
$$N \cap (N^\al)' = \Jte (N_2 \cap N') \Jte = \C.$$
\end{proof}

As a complement to theorem~\ref{52} we prove the following easy
result. The terminology is taken from \cite{NT}.
\begin{proposition} \label{semidual}
Let $\al$ be an action of \qu on $N$. Then we call $\al$ semidual
when there exists a unitary $v \in B(H) \ot N$ satisfying $(\io
\ot \al)(v) = v_{13}V_{12}^*$.
\begin{itemize}
\item Every dual action is semidual.
\item Every semidual action is integrable and the $*$-homomorphism
$\rho$ from theorem~\ref{52} is faithful.
\end{itemize}
\end{proposition}
\begin{proof}
Let us first prove the first statement.
Denote with $\alh$ the dual action, which is an action of
$(\Mh,\dehop)$ on $\cros$. Because $\vfih$ is the right Haar
weight of $(\Mh,\dehop)$, the role of $V$ is played by $\Si
\hat{W}^* \Si = W$. So we have to find a unitary $v \in B(H) \ot
(\cros)$ satisfying $(\io \ot \alh)(v) = v_{13} W_{12}^*$. Then it
is clear that we can take $v= W^* \ot 1$ and so $\alh$ is
semidual.

To prove the second part suppose that $v \in B(H) \ot N$ is
unitary and $(\io \ot \al)(v) = v_{13} V_{12}^*$. Define the
isomorphism $\Psi : B(H) \ot N \recht B(H) \ot N$ by
$\Psi(z)=vzv^*$. Using the notation of theorem~\ref{26} we get
that $\mu(\Psi(z)) = (\io \ot \Psi)\ga(z)$ for all $z \in B(H) \ot
N$. So the action $\mu$ of \qu on $B(H) \ot N$ is isomorphic with
the action $\ga$, which is integrable because it is isomorphic
with the bidual action $\alhh$. Hence $\mu$ is integrable, and so
$\al$ is integrable.

Fix now a \nsf weight $\te$ on $N$ and represent $N$ on the
GNS-space of $\te$ such that $(K,\io,\late)$ is a
GNS-construction. Let $N_2 = \Jte (N^\al)' \Jte$ be the basic
construction from $N^\al \subset N$ and let $\rho : \cros \recht
N_2$ be the $*$-homomorphism from theorem \ref{52}.
Then define $w=(\Jh \ot \Jte)v(\Jh \ot \Jte)$ and define
$$\eta : N_2 \rightarrow B(H \ot K) : \eta(z) = U w^* (1 \ot z) w
U^* \tekst{for all} z \in N_2.$$
Because $w \in B(H) \ot N'$ we have
$$\eta(x) = U (1 \ot x) U^* = \al(x)$$
for all $x \in N$. Further we have $(\io \ot \al)(v) =
v_{13}V_{12}^*$ and so $U_{23} v_{13} U_{23}^* = v_{13} V_{12}^*$.
Putting $\Jh \ot \Jh \ot \Jte$ around this equation and using that
$V = (\Jh \ot \Jh) \Si W^* \Si (\Jh \ot \Jh)$ (see \cite[2.15]{KV3}), we get
$$U_{23}^* w_{13} U_{23} = w_{13} (\Si W \Si)_{12}.$$
Flipping the first two legs of this equation and rewriting it
we get
$$w_{23}^* U_{13}^* w_{23} = W_{12} U_{13}^*.$$
From this it follows that
$$U_{23} w_{23}^* U_{13}^* w_{23} U_{23}^* = U_{23} W_{12}
U_{13}^* U_{23}^* = U_{23} W_{12} (\de \ot \io)(U^*) = U_{23}
U_{23}^* W_{12} = W_{12}.$$
Then we get for all $\om \in M_*$ that
$$\eta \bigl( (\om \ot \io)(U^*) \bigr) = (\om \ot \io \ot \io)(U_{23} w_{23}^* U_{13}^* w_{23}
U_{23}^*) = (\om \ot \io)(W) \ot 1.$$
Hence we may conclude that $\eta \na \rho = \io$ and so $\rho$ is
faithful.
\end{proof}

\sectie{Minimal actions and outer actions} \label{minimal}
In definition \ref{outer} we already defined the notion of an outer action.
In the literature one usually encounters the notion of outer
action when dealing with discrete group actions and one encounters
the notion of minimal action when dealing with compact group
actions. In this section we will prove how both notions can be
linked in a locally compact quantum group setting. We will also
prove a generalization of the main theorem of Yamanouchi,
\cite{Ya}: when working on separable Hilbert spaces,
we prove that every integrable outer action with infinite fixed
point algebra is a dual action.

The following definition appears in \cite[4.3]{ILP} when
dealing with actions of compact Kac algebras.

\begin{definition}
An action $\al$ of \qu on $N$ is called minimal when
$$N \cap (N^\al)' = \C \quad\tekst{and}\quad \{(\io \ot \om)\al(x)
\mid \om \in N_*, x \in N \}^{\prime\prime}=M.$$
\end{definition}

We will prove the following result.
\begin{proposition} \label{outermin}
Let $\al$ be an action of \qu on $N$.
\begin{itemize}
\item If $\al$ is minimal, then $\al$ is outer.
\item If $\al$ is outer and integrable, then $\al$ is minimal.
\end{itemize}
\end{proposition}
\begin{proof}
Let $\al$ be minimal. Let $z \in (\cros) \cap \al(N)'$. Then
certainly $z \in (B(H) \ot N) \cap (\C \ot N^\al)'$ and hence $z \in
B(H) \ot \C$ by minimality. We now claim that for $x \in B(H)$ we
have $x \ot 1 \in \cros$ if and only if $x \in \Mh$. Suppose $x
\ot 1 \in \cros$. It is clear that for every $z \in \cros$ we have
$(\io \ot \al)(z)= V_{12}z_{13}V_{12}^*$. So we get $(x \ot 1) V =
V (x \ot 1)$. From this it follows that $x \in \Mh$. So we may
conclude that $z = x \ot 1$, where $x \in \Mh$. Because $z \in
\al(N)'$ we get that $(x \ot 1) \al(y) = \al(y) (x \ot 1)$ for all
$y \in N$. By minimality we get $x \in M'$. But then $x \in M' \cap \hat{M} = \C$
and so $z \in \C$. Hence $\al$ is outer.

Let now $\al$ be outer and integrable. Choose a \nsf weight $\te$ on $N$ and
represent $N$ on the GNS-space of $\te$. Let $\Jte$ denote the
modular conjugation of $\te$ and let $N_2 = \Jte (N^\al)' \Jte$ be
the basic construction from $N^\al \subset N$.
Let $\rho$ be the $*$-homomorphism given in theorem \ref{52}. Then $\rho$ is faithful
because $\cros$ is a factor. Because $\rho$ is
an isomorphism we get $N_2 \cap N' = \C$ and so
$$N \cap (N^\al)' = \Jte (N_2 \cap N') \Jte = \C.$$
Next we claim that $\bigl( \al(N) \cup \C \ot N' \bigr)\dpr= M \ot B(K)$.
Because, by theorem \ref{26}, $B(H) \ot N = \bigl( \cros \cup M \ot \C \bigr)\dpr$,
we get
\begin{align*}
B(H) \ot \bigl( \al(N) \cup \C \ot N' \bigr)\dpr &= \bigl( (\io \ot \al)(B(H) \ot
N) \cup \C \ot \C \ot N' \bigr)\dpr \\
&= \bigl( (\io \ot \al)(\cros) \cup M \ot \C \ot N' \bigr)\dpr \\
&= V_{12} \bigl( (\cros)_{13} \cup V^*(M \ot \C)V \ot N' \bigr)\dpr
V_{12}^*.
\end{align*}
When $\Jtil$ denotes the modular conjugation of the dual weight
$\tetil$, we already observed in the proof of proposition~\ref{depth} that
$B(H) \ot N = \Jtil \al(N)' \Jtil$. Then the outerness of $\al$
implies that $B(H) \ot N \cap (\cros)' = \C$ and so
$$\bigl( \C \ot N' \cup \cros \bigr)\dpr = B(H) \ot B(K).$$
Then we may conclude from the previous computation that
\begin{align*}
B(H) \ot \bigl( \al(N) \cup \C \ot N' \bigr)\dpr &= V_{12} \bigl( B(H) \ot \C \ot
B(K) \cup V^*(M \ot \C)V \ot \C \bigr)\dpr V_{12}^* \\
&= \bigl( V(B(H) \ot \C)V^* \ot B(K) \cup M \ot \C \ot \C \bigr)\dpr \\
&=\bigl( \de(M) \ot B(K) \cup (\Mh \cup M) \ot \C \ot \C \bigr)\dpr \\
&=B(H) \ot M \ot B(K)
\end{align*}
where we have used that $V \in \Mh' \ot M$, $(\Mh \cup M)\dpr =
B(H)$ and $\bigl( \de(M) \cup B(H) \ot \C \bigr)\dpr = B(H) \ot M$. Then our
claim follows and hence it is clear
that
$$\{(\io \ot \om)\al(x)
\mid \om \in N_*, x \in N \}^{\prime\prime}=M.$$
So $\al$ is minimal.
\end{proof}
We will now give an example of an outer action which is not
minimal.
\begin{counterex} \label{example}
There exists an action $\al$ of $\Z$ on a $II_1$-factor $N$ such
that $\al$ is outer and $N^\al = \C$. Then $\al$ is clearly not
minimal, and neither can $\C \ot N^\al \subset \al(N) \subset
\cros$ be the basic construction.
\end{counterex}
\begin{proof}
Let $G$ be the free group with a countably infinite number of
generators $\{a_n \mid n \in \Z\}$. It is well known that the free group factor
$N=\cL(G)$ is a $II_1$-factor. Let $\be$ be the automorphism of
$G$ satisfying $\be(a_n)=a_{n+1}$ for all $n \in \Z$. Let $\al$ be
the automorphism of $N$ satisfying $\al(\sla_g) = \sla_{\be(g)}$
for all $g \in G$. Define the automorphism group $(\al_n)_{n \in
\Z}$ in the usual way by $\al_n = \al^n$ for all $n \in \Z$. It is
easy to verify that $\al$ is a free action and hence $\al$ is outer
(see \cite[def. 1.4.2 and prop. 1.4.4]{JS}). Further it is easy to
check that $N^\al = \C$.
\end{proof}
We conclude this section with a generalization of the main theorem
of Yamanouchi \cite{Ya}. It is remarkable that the proof of our result is
much more easy then Yamanouchi's proof. In \cite{Ya} the following
result is proved for minimal actions of compact Kac algebras,
which are automatically integrable because the Haar weight is
finite.
\begin{proposition}
Let $\al$ be an action of \qu on $N$. Suppose that both $M$ and
$N$ are $\si$-finite von Neumann algebras (i.e. with separable
preduals). If the action $\al$ is minimal and integrable and if
$N^\al$ is infinite, then $\al$ is a dual action.
\end{proposition}
\begin{proof}
Consider the action $\be$ of \qu on $\Ntil=B(H) \ot N \ot
M_2(\C)$ given by
$$\be \begin{pmatrix} x & y \\ z & r \end{pmatrix} = \begin{pmatrix} \mu(x) & \mu(y) \Vtil^*
\\ \Vtil \mu(z) & \ga(r) \end{pmatrix},$$
for $x,y,z,r \in B(H) \ot N$.
Here we used the notations of theorem~\ref{26}: $\mu(x) = (\si \ot
\io)(\io \ot \al)(x)$, $\ga(r) = \Vtil \mu(r) \Vtil^*$ and $\Vtil
= \Si V^* \Si \ot 1$. Let us define now
$$\cJ = \{ x \in B(H) \ot N \mid (\io \ot \al)(x) = x_{13}
V_{12}^*\}.$$
Using matrix notation and referring to theorem \ref{26} and \ref{27}, it is then clear that $x \in \Ntil^\be$ if
and only if $x_{11} \in B(H) \ot N^\al$, $x_{22} \in \cros$ and
$x_{12},x_{21}^* \in \cJ$.

Choose a \nsf weight $\te$ on $N$ and represent $N$ on the
GNS-space of $\te$ such that $(K,\io,\late)$ is a
GNS-construction.
Then we fix $z \in \cN_{(\psi \ot \io)\al}$ and $\xi \in H$ and we claim that the
element $x \in B(H \ot K)$ defined by
$$x:= (\Ga \ot \io)\al(z) (\te_\xi^* \ot 1)$$
belongs to $\cJ^*$. Here we used the notation $\Ga \ot \io$
introduced in the proof of theorem~\ref{52}.
To prove our claim we observe that for all $b \in \cN_\psi$, $y \in \Nte$
and $\eta \in H$
$$(\te_{\Ga(b)}^* \ot 1) x (\eta \ot \late(y)) = \langle \eta,\xi
\rangle (\te_{\Ga(b)}^* \ot 1)(\Ga \ot \late)(\al(z)(1 \ot y)) = \langle \eta,\xi
\rangle (\psi \ot \io) \bigl( (b^* \ot 1)\al(z) \bigr) \late(y).$$
We can conclude that
$$(\om_{\eta,\Ga(b)} \ot \io)(x) = \langle \eta,\xi
\rangle  (\psi \ot \io) \bigl( (b^* \ot 1) \al(z) \bigr).$$
So $x \in B(H) \ot N$ and for all $\eta \in H, b \in \cN_\psi$ and
$\om \in N_*$ we have
\begin{align*}
(\om_{\eta,\Ga(b)} \ot \io \ot \om)(\io \ot \al)(x) &= (\io \ot
\om) \al \bigl( \langle \eta,\xi
\rangle  (\psi \ot \io) \bigl( (b^* \ot 1) \al(z) \bigr) \bigr) \\
&=\langle \eta,\xi \rangle (\psi \ot \io) \bigl( (b^* \ot 1) \de((\io \ot
\om)\al(z)) \bigr) \\
&=\langle \eta,\xi \rangle (\om_{\Ga ( (\io \ot \om)\al(z)),\Ga(b)}
\ot \io)(V).
\end{align*}
Next we observe that for all $y \in \cN_\psi$
\begin{align*}
\langle \langle \eta,\xi \rangle \Ga \bigl( (\io \ot \om) \al(z) \bigr), \Ga(y)
\rangle &= \langle \eta,\xi \rangle \om (\psi \ot \io) \bigl( (y^* \ot
1) \al(z) \bigr) \\
&= \om \bigl( (\om_{\eta,\Ga(y)} \ot \io)(x) \bigr) = \langle (\io \ot \om)(x)
\eta, \Ga(y) \rangle.
\end{align*}
Inserting this in the computation above we get that
$$(\om_{\eta,\Ga(b)} \ot \io \ot \om)(\io \ot \al)(x) = (\om_{(\io
\ot \om)(x) \eta, \Ga(b)} \ot \io)(V) = (\om_{\eta,\Ga(b)} \ot \io \ot
\om)(V_{12} x_{13}).$$
Then it follows that $x \in \cJ^*$.

So we see that $\cJ \neq \{0\}$. Because $\al$ is minimal we also
have that $\al$ is outer by proposition~\ref{outermin}. In particular $\cros$ is a
factor. Also $N^\al$ is a factor. Because $\cJ \neq \{0\}$ we then get immediately that
$\Ntil^\be$ is a factor. Because $M$ is supposed to be $\si$-finite,
the Hilbert space $H$ is separable. So $\Ntil^\be$ is $\si$-finite.
Denoting with $e_{ij}$ the matrix units in $M_2(\C)$ we see that
the projections $1 \ot e_{11}$ and $1 \ot e_{22}$ both belong to
$\Ntil^\be$. Because $\C \ot N^\al \subset \cros$ both projections
are infinite. Hence they are equivalent in the $\si$-finite factor
$\Ntil^\be$. Take $w \in \Ntil^\be$ such that $w^*w = 1 \ot
e_{22}$ and $w w^* = 1 \ot e_{11}$. Then there exists a unitary $v
\in \cJ$ such that $w = v \ot e_{12}$.

Now we can consider the isomorphism
$$\Psi: B(H) \ot N \recht B(H) \ot N : \Psi(z)=v^* z v.$$
It is easy to check that $(\io \ot \Psi)\mu(z)=\ga(\Psi(z))$ for
all $z \in B(H) \ot N$. So the actions $\mu$ and $\ga$ are
isomorphic. Because $\ga$ is isomorphic to the bidual action
$\alhh$ by theorem \ref{26}, we get that $\mu$ is a dual action. Because
$N^\al$ is properly infinite and because $H$ is a separable Hilbert space
we get that the action $\al$ on $N$ is isomorphic with the action
$\mu$ on $B(H) \ot N$. So $\al$ is a dual action.
\end{proof}

\sectie{Appendix}
In this appendix we collect four technical results which do not
have anything to do with actions. The first three results are  general
results on locally compact quantum groups and the last one deals
with \nsf weights on a \vna. We will use freely the notations
introduced in the introduction.

\begin{proposition} \label{a1}
Let \qu be a locally compact quantum group. For every $\xi \in H$
and $b \in \cT_\vfi$ we have
$$\pih(\om_{\xi,\la(b)}) \in \Nfih \tekst{and}
\lah(\pih(\om_{\xi,\la(b)})) = J \si_{i/2}(b) J \xi.$$
Moreover $\lspan \{ \pih(\om_{\xi,\la(b)}) \mid \xi \in H, b \in
\cT_\vfi \}$ is a \strong--norm core for $\lah$.
\end{proposition}
\begin{proof}
The first statement follows easily from the definition of $\vfih$. Let $x \in \Nfi$, then
$$\om_{\xi,\la(b)}(x^*)=\langle x^* \xi,\la(b) \rangle = \langle
J \si_{i/2}(b) J \xi, \la(x) \rangle.$$
So we get the first statement. To prove the second one we define
$$\cL = \{ a \in \Nfi \mid \; \text{there exists}\quad \om \in M_*
\tekst{such that} \om(x)=\vfi(xa) \tekst{for all} x \in \Nfi^*
\}.$$
It is clear that for $a \in \cL$ such a $\om \in M_*$ is
necessarily unique. We denote it with $a \vfi$. Then for every $a
\in \cL$ we have $\pih(a \vfi) \in \Nfih$ and $\lah(\pih(a
\vfi))=\la(a)$. Define $\cD_0 = \{\pih(a\vfi) \mid a \in \cL\}$.
We claim that $\cD_0$ is a \strong--norm core for $\lah$.

Denote with $\cD$ the domain of the \strong--norm closure of the
restriction of $\lah$ to $\cD_0$.

Let $a \in \cL$ and $t \in \R$. Define $b = \tau_t(a) \sde^{-it}$.
Then $b \in \Nfi$ and for all $x \in \Nfi^*$ we have
$$\vfi(xb)=\vfi(x\tau_t(a) \sde^{-it}) =\nu^t \vfi(\sde^{-it} x
\tau_t(a)) = \vfi(\sde^{-it} \tau_{-t}(x) a) = (a\vfi)(\sde^{-it}
\tau_{-t}(x)) = (\rho_t(a\vfi))(x)$$
where we used the notation of \cite[8.7]{KV2}. So $b \in \cL$ and $b \vfi =
\rho_t(a \vfi)$. Hence
$$\sih_t(\pih(a\vfi)) = \pih(\rho_t(a\vfi)) = \pih(b \vfi) \in
\cD_0.$$
So we get that $\cD_0$ is invariant under $\sih_t$. Then it is
easy to conclude that $\cD$ is invariant under $\sih_t$ for all $t
\in \R$.

Let now $\om \in M_*$ and suppose that there exists a $\mu \in
M_*$ such that $\mu(x) = \om(S^{-1}(x))$ for all $x \in
\cD(S^{-1})$. Let $a \in \cL$. Define $b = (\mu \ot \io)\de(a)$.
Then $b \in \Nfi$ and for all $x \in \Nfi^*$ we have
\begin{align*}
\vfi(xb) &= \vfi \bigl( x (\mu \ot \io)\de(a) \bigr) = \mu \bigl( (\io \ot \vfi) \bigl( (1
\ot x) \de(a) \bigr) \bigr) = \om \bigl( (\io \ot \vfi) \bigl( \de(x)(1 \ot a) \bigr) \bigr) \\
&= \vfi \bigl( (\om \ot \io)\de(x) a \bigr) = (\om \ot a\vfi)\de(x).
\end{align*}
So we see that $b \in \cL$ and $b \vfi = (\om \ot a \vfi) \de$.
Then we may conclude that
$$\pih(\om) \pih(a \vfi) = \pih(b \vfi) \in \cD_0.$$
Because such elements $\pih(\om)$ form a \strong dense subset of
$\Mh$ it is easy to conclude that $\cD$ is a left ideal in $\Mh$.

Because $\cD$ is a \strong dense left ideal of $\Mh$, invariant under $\sih$
and because $\cD \subset \Nfih$, we may conclude that $\cD$ is a
\strong--norm core for $\lah$. But then $\cD = \Nfih$ and we have
proven our claim.

Then it follows easily that also
$$\lspan \{ \pih(ab \vfi) \mid a \in \cL, b \in \cT_\vfi \}$$
is a \strong--norm core for $\lah$. This last space equals
$$\lspan \{ \pih(\om_{\la(a),\la(b)}) \mid a \in \cL, b \in
\cT_\vfi \}$$
and so the proposition is proven.
\end{proof}
For completeness we also include the following easy result.
\begin{proposition} \label{a2}
Let \qu be a locally compact quantum group. For every $a \in
\Nfi$, $\xi \in \cD(\sde^{1/2})$ and $\eta \in H$ we have
$(\io \ot \om_{\xi,\eta})\de(a) \in \Nfi$ and
$$\la \bigl( (\io \ot \om_{\xi,\eta})\de(a) \bigr) = (\io \ot \om_{\sde^{1/2}
\xi,\eta})(V) \la(a).$$
\end{proposition}
\begin{proof}
Let $(e_n)$ be the sequence of operators defined in the proof of \cite[7.6]{KV2}. Because
$\de(\sde)=\sde \ot \sde$ it is clear that
$$\bigl( (\io \ot \om_{\xi,\eta})\de(ae_n) \bigr) \sde^{-1/2} \subset (\io \ot
\om_{\sde^{1/2} \xi,\eta}) \de \bigl( a (\sde^{-1/2} e_n) \bigr).$$
Because $a (\sde^{-1/2} e_n) \in \Npsi$ we have $(\io \ot
\om_{\sde^{1/2} \xi,\eta}) \de \bigl( a (\sde^{-1/2} e_n) \bigr) \in \Npsi$. We
know that $\vfi=\psi_{\sde^{-1}}$, so that $(\io \ot
\om_{\xi,\eta})\de(ae_n) \in \Nfi$ and
\begin{align*}
\la \bigl( (\io \ot \om_{\xi,\eta})\de(ae_n) \bigr) &= \Ga \bigl( (\io \ot
\om_{\sde^{1/2} \xi,\eta}) \de \bigl( a (\sde^{-1/2} e_n) \bigr) \bigr) \\
&=(\io \ot \om_{\sde^{1/2} \xi,\eta})(V) \Ga \bigl( a(\sde^{-1/2}e_n) \bigr) \\
&=(\io \ot \om_{\sde^{1/2} \xi,\eta})(V) \la(ae_n).
\end{align*}
Because $\la$ is \strong--norm closed, the conclusion follows.
\end{proof}

We also need the following technical result.
\begin{proposition} \label{extra}
Let \qu be a locally compact quantum group and let $x \in M$.
Suppose that there exists a vector $\eta \in H$ such that
$$\om(x^*) = \langle \xi(\om) , \eta \rangle$$
for all $\om \in \cI$. Then $x \in \Nfi$ and $\la(x)=\eta$.
\end{proposition}
\begin{proof}
Let $\om \in \cI$ and $y \in \cT_\vfi$. Then we have for all $a
\in \Nfi$ that
$$(\om y)(a^*) = \om \bigl( (ay^*)^* \bigr) = \langle \xi(\om),\la(ay^*)
\rangle = \langle J \si_{i/2}(y^*) J \xi(\om) , \la(a) \rangle.$$
So we get that $\om y \in \cI$ and $\xi(\om y) = J \si_{i/2}(y^*)
J \xi(\om)$. Take now a net $(e_\al)$ in $\cT_\vfi$ such that
$\si_z(e_\al) \recht 1$ in the \strong topology for all $z \in
\C$. Then we have for all $\om \in \cI$
$$\langle \xi(\om),\la(x e_\al) \rangle = \om(e_\al^* x^*) = (\om
e_\al^*)(x^*) = \langle \xi(\om e_\al^*), \eta \rangle = \langle
\xi(\om), J \si_{i/2}(e_\al)^* J \eta \rangle.$$
Hence $\la(x e_\al) = J \si_{i/2}(e_\al)^* J \eta$ for all $\al$.
Because $\la$ is \strong--norm closed we get $x \in \Nfi$ and
$\la(x)=\eta$.
\end{proof}

The following result is probably well known, but we could not find
it in the literature.
\begin{proposition} \label{a3}
Let $\te$ be a \nsf weight on a \vna $N$ with GNS-construction
$(H,\pi,\la)$. Suppose that
\begin{itemize}
\item $\cD$ is a weakly dense left ideal in $N$ with $\cD \subset
\Nte$.
\item $K$ is a Hilbert space and $\la_0: \cD \recht K$ is a linear
map such that $\la_0(\cD)$ is dense in $K$.
\item $\pi_0$ is normal representation of $N$ on $K$ such that
$\pi_0(x)\la_0(y) = \la_0(xy)$ for all $x \in N$ and $y \in \cD$.
\item $\cV$ is an isometry from $K$ to $H$ such that $\cV
\la_0(x)=\la(x)$ for all $x \in \cD$.
\item $\la_0$ is \strong--norm closed.
\end{itemize}
Then there exists a unique \nsf weight $\mu$ on $N$ such that
$\cN_\mu = \cD$ and $(K,\pi_0,\la_0)$ is a GNS-construction for
$\mu$. In particular $\mu$ is a restriction of $\te$, which means
that for every $x \in \cM_\mu^+$ we have $x \in \cM_\te^+$
and $\mu(x)=\te(x)$.
\end{proposition}
\begin{proof}
Because $\cV$ is an isometry, $\la_0$ is injective. Define $\cU =
\la_0(\cD \cap \cD^*)$. Then $\cU$ is a dense subspace of $K$. We
make $\cU$ into a $*$-algebra by using $\la_0$ and the $*$-algebra
structure on $\cD \cap \cD^*$. We claim that $\cU$ is a left
Hilbert algebra. The only non-trivial point is to prove that the
map $\la_0(x) \mapsto \la_0(x^*)$ for $x \in \cD \cap \cD^*$ is
closable. But, suppose that $(x_n)$ is a sequence in $\cD \cap
\cD^*$ such that $\la_0(x_n) \recht 0$ and $\la_0(x_n^*) \recht
\xi \in K$. Applying $\cV$ we get $\la(x_n) \recht 0$ and $\la(x_n^*) \recht \cV
\xi$. Because $\te$ is a \nsf weight we get that $\cV \xi = 0$ and
so $\xi = 0$. This gives our claim.

It is clear that the \vna generated by the left Hilbert algebra
$\cU$ is $\pi_0(N)$. Because $\la_0$ is injective we have that
$\pi_0$ is injective. So the \nsf weight on $\pi_0(N)$ which is canonically
associated with $\cU$ can be composed with $\pi_0$ to obtain a
\nsf weight $\mu$ on $N$. Let $(K,\pi_0,\la_1)$ be the canonically
associated GNS-construction.

Then, by definition of $\mu$, every $x \in \cD \cap \cD^*$ will
belong to $\cN_\mu$ and $\la_1(x)=\la_0(x)$. Let now $x \in \cD$.
Take a net $(e_\al)$ in $\cD$ such that $e_\al \recht 1$ \strong.
Then we have $e_\al^* x \recht x$ \strong, $e_\al^* x \in \cD \cap
\cD^*$ and
$$\la_1(e_\al^* x) = \la_0(e_\al^* x) = \pi_0(e_\al^*) \la_0(x)
\recht \la_0(x).$$
Because $\la_1$ is \strong--norm closed we get $x \in \cN_\mu$ and
$\la_1(x)=\la_0(x)$.

Conversely, suppose $x \in \cN_\mu$. By the main theorem of \cite{H3}
there exists a net $(x_\al)$ in $\cD \cap \cD^*$ such
that $\|x_\al\| \leq \|x\|$ for all $\al$, $x_\al \recht x$
\strong and $\la_1(x_\al) \recht \la_1(x)$ in norm. But
$\la_1(x_\al)=\la_0(x_\al)$ for every $\al$. Because $\la_0$ is
\strong--norm closed we get that $x \in \cD$. This concludes our
proof.
\end{proof}

\end{document}